\documentclass[openright,11pt,two0side,a4paper,leqno]{article}
\usepackage[latin1]{inputenc} 
\usepackage{amsmath,amssymb} 
\usepackage{graphicx}
\usepackage{theorem}

\newlength{\larg}
\newcommand{\pv}{\noindent \emph{Proof. }}
\newcommand{\cqfd}{\hfill $\Box$}
\newcommand{\fini}{\hfill $\diamond$}
\newenvironment{demo}{\pv}{\cqfd \\ }

\newcommand{\ie}{\emph{i.e. }}
\newcommand{\eg}{\emph{e.g. }}
\newcommand{\resp}{\emph{resp. }}

\newcommand{\barre}{\overline{ ^{\ \: }} \ }

\newcommand{\set}[1]{ \{ #1 \}}

\newcommand{\mC}{{\mathbb C}}
\newcommand{\mN}{{\mathbb N}}

\newcommand{\mR}{{\mathbb R}}
\newcommand{\mZ}{{\mathbb Z}}

\newcommand{\ds}{\displaystyle}
\newcommand{\interv}[2]{[\! [ #1  ;  #2 ] \!]}

\newcommand{\hsl}[1]{{\widehat{\mathfrak{sl}}}_{#1}}
\newcommand{\Uq}{U_q(\hsl{n})}

\newcommand{\Fq}[1][\bd{s}_l]{\boldsymbol{\operatorname{F}}_q[\textbf{s}_l]}

\newcommand{\Sr}{\mathfrak{S}_r}

\newcommand{\lambdal}{\boldsymbol{\lambda}_l}
\newcommand{\mul}{\boldsymbol{\mu}_l}
\newcommand{\nul}{\boldsymbol{\nu}_l}
\newcommand{\esselle}{\textbf{s}_l}
\newcommand{\Pilm}{\Pi^l_m}

\newcommand{\cont}{\operatorname{cont}}
\newcommand{\hd}{\operatorname{hd}}
\newcommand{\haut}{\operatorname{ht}}
\newcommand{\res}{\operatorname{res}}
\newcommand{\tl}{\operatorname{tl}}

\newtheorem{thm}{Theorem}[section]
\newtheorem{lemma}[thm]{Lemma}
\newtheorem{prop}[thm]{Proposition}
\newtheorem{cor}[thm]{Corollary}
\newtheorem{conj}[thm]{Conjecture}

{\theorembodyfont{\rmfamily} \newtheorem{remark}[thm]{Remark}}   
{\theorembodyfont{\rmfamily} \newtheorem{example}[thm]{Example}} 
{\theorembodyfont{\rmfamily} \newtheorem{definition}[thm]{Definition}}
{\theorembodyfont{\rmfamily} \newtheorem{notation}[thm]{Notation}}

\title{A conjecture for $q$-decomposition matrices of cyclotomic $v$-Schur algebras}
\author{Xavier YVONNE}

\begin{document}

\maketitle

\begin{abstract} The Jantzen sum formula for cyclotomic $v$-Schur algebras yields an identity for some $q$-analogues of the decomposition matrices 
of these algebras. We prove a similar identity for matrices of canonical bases of higher-level Fock spaces. We conjecture then that those 
matrices are actually identical for a suitable choice of parameters. In particular, we conjecture that decomposition matrices of cyclotomic 
$v$-Schur algebras are obtained by specializing at $q=1$ some transition matrices between the standard basis and the canonical basis of a Fock space. 
\end{abstract}

\section{Introduction}
In order to study representations of the Ariki-Koike algebra associated to the complex reflection group $G(l,1,m)$, Dipper, James and Mathas 
introduced in 1998 the cyclotomic $v$-Schur algebra \cite{DJM}. This algebra depends on the two integers $l$ and $m$ and on some deformation 
parameters $v,u_1,\ldots,u_l$. When $l=1$, the cyclotomic $v$-Schur algebra coincides with the $v$-Schur algebra of \cite{DJ}. It is an open problem 
to calculate the decomposition matrix of a cyclotomic $v$-Schur algebra whose parameters are powers of a given $n$-th root of unity. To this aim, 
James and Mathas proved, for cyclotomic $v$-Schur algebras, an important formula: the Jantzen sum formula \cite{JM}. Given a Jantzen filtration for 
Weyl modules, one can define a $q$-analogue $D(q)$ of the decomposition matrix; the coefficients of $D(q)$ are graded decomposition numbers
of the composition factors of Weyl modules (see Definition \ref{definition_Dq}). The Jantzen sum formula is equivalent to the identity 
$D'(1)=J^{\lhd}D(1)$, where $J^{\lhd}$ is a matrix of $\wp$-adic valuations of factors of some Gram determinants (see Theorem 
\ref{thm_Jantzen_formula} and Corollary \ref{cor_Jantzen_formula}). \\  
\indent Let $\Delta(q)$ be the matrix of the canonical basis of the degree $m$ homogeneous component of a Fock representation of level $l$ of $\Uq$ 
\cite{U2}. Uglov provided in \cite{U2} an algorithm for computing $\Delta(q)$. \\
\indent In view of Ariki's theorem for Ariki-Koike algebras \cite{A2}, it seems natural to conjecture that for a suitable choice of parameters, 
one has $D(q)=\Delta(q)$. This would provide an algorithm for computing decomposition matrices of cyclotomic $v$-Schur algebras. Varagnolo and 
Vasserot \cite{VV} proved for $l=1$ that $D(1)=\Delta(1)$. Moreover, Ryom-Hansen showed that this conjecture (still for $l=1$) is compatible with 
the Jantzen-Schaper formula \cite{Ry}. Passing to higher level $l \geq 1$ requires the introduction of an extra parameter 
$\esselle = (s_1,\ldots,s_l) \in \mZ^l$, called \emph{multi-charge}; this $l$-tuple parametrizes the Fock space of level $l$ introduced by Uglov. 
We say that $\esselle$ is \emph{$m$-dominant} if for all $1 \leq d \leq l-1$, we have $s_{d+1}-s_d \geq m$. In this case, we
conjecture that $D(q)=\Delta(q)$. Here, $D(q)$ comes from a Jantzen filtration of the Weyl modules of the cyclotomic $v$-Schur algebra
$\mathcal{S}_{\mC}=\mathcal{S}_{\mC,m}(\zeta;\zeta^{s_1},\ldots,\zeta^{s_l})$ with $\zeta:=\exp(\frac{2i\pi}{n})$. Note that for any choice of roots
of unity $\zeta^{\mathtt{r}_1},\ldots,\zeta^{\mathtt{r}_l}$ (that is, for any $\mathtt{r}_1,\ldots,\mathtt{r}_l \in \mZ/n\mZ$) and any $m$ we can 
find an $m$-dominant multi-charge $\esselle=(s_1,\ldots,s_l)$ such that $\zeta^{s_d}=\zeta^{\mathtt{r}_d}$ ($1 \leq d \leq l$). Therefore, putting 
$q=1$, our conjecture gives an algorithm for calculating the decomposition matrix of an arbitrary cyclotomic $v$-Schur algebra 
$\mathcal{S}_{\mC}=\mathcal{S}_{\mC,m}(\zeta;\zeta^{s_1},\ldots,\zeta^{s_l})$. Such a conjecture is new even for type $B_m$ (case $l=2$). \\  
\indent Our conjecture is supported by the following theorem. We define in a combinatorial way a matrix $J^{\prec}$ for any multi-charge $\esselle$;
if $\esselle$ is $m$-dominant, then our matrix $J^{\prec}$ coincides with the matrix $J^{\lhd}$ of the Jantzen sum formula.
We show then that for any multi-charge $\esselle$, we have $\Delta'(1)=J^{\prec} \Delta(1)$ (Theorem \ref{thm_Delta_J}). \\
\indent The proof of our theorem relies on a combinatorial expression for the derivative at $q=1$ of the matrix $A(q)$, where $A(q)$ is the matrix
of the Fock space involution used for defining $\Delta(q)$. Namely, we show that $A'(1)=2J^{\prec}$ (Theorem \ref{thm_A_J}). The 
coefficients of $A(q)$ are some analogues for Fock spaces of Kazhdan-Lusztig $R$-polynomials $R_{x,y}(q)$ for Hecke algebras. The classical
computation of $R'_{x,y}(1)$ was made in \cite{GJ}, in relation with the Kazhdan-Lusztig conjecture for multiplicities of composition factors of 
Verma modules. \\ 

\noindent \emph{Acknowledgments.} I would like to thank Nicolas Jacon and my advisor Bernard Leclerc for inspiring discussions about Ariki-Koike
algebras. I also would like to thank Bernard Leclerc for his assistance and constant advice when I was writing this article. At last, I thank
Andrew Mathas, Hyohe Miyachi and the referee for their comments.

\begin{notation} \label{notation_S_r}
Let $\mN$ (\resp $\mN^{*}$) denote the set of nonnegative (\resp positive) integers, and for $a,\,b \in \mR$ denote by $\interv{a}{b}$ the discrete 
interval $[a  ,  b] \cap \mZ$. Throughout this article, we fix three integers $n,\,l,\,m \geq 1$. Let $\Pi$ be the set of partitions of any integer
and $\Pilm$ be the set of $l$-multi-partitions of $m$. The Coxeter group of type $A_{r-1}$ (with $r \in \mN^*$) is the symmetric group
$\Sr = \langle \sigma_i=(i,i+1) \mid 1 \leq i \leq r-1 \rangle$. Let $\ell$ be the length function on $\Sr$ and $\omega$ be the unique element of 
maximal length in $\Sr$. \fini
\end{notation}

\vspace{5mm}
\begin{center}
\scshape \Large PART A: Statement of results 
\end{center}
\vspace{1mm}

\section{Statement of results} 
\subsection{The Jantzen sum formula} 

\begin{definition}[\cite{AK,BM}]
Let $R$ be a principal ideal domain. Let $v$ be an invertible element of $R$ and $u_1,\ldots,u_l \in R$. The 
\emph{Ariki-Koike algebra}, denoted by 
\begin{equation}
\mathcal{H}=\mathcal{H}_R=\mathcal{H}_{R,m}(v;u_1,\ldots,u_l),
\end{equation}
is the algebra defined over $R$ with generators $T_0,\ldots,T_{m-1}$ and relations
\begin{equation}
\left\{
\begin{array}{rclc}    
(T_0-u_1) \cdots (T_0-u_l) &=& 0, \\
T_0 T_1 T_0 T_1 &=& T_1 T_0 T_1 T_0, \\
(T_i+1)(T_i-v) &=& 0 & (1 \leq i \leq m-1), \\
T_i T_{i+1} T_i &=& T_{i+1} T_i T_{i+1} & (1 \leq i \leq m-2), \\
T_i T_j &=& T_j T_i & (0 \leq i < j-1 \leq m-2).
\end{array} \right.
\end{equation}
\fini
\end{definition}

\! \,Following \cite{DJM}, let 
\begin{equation}
\mathcal{S}=\mathcal{S}_R=\mathcal{S}_{R,m}(v;u_1,\ldots,u_l)
\end{equation}
be the cyclotomic $v$-Schur algebra associated to $\mathcal{H}$. Dipper, James and Mathas (see \cite[Theorem 6.12]{DJM}) showed that $\mathcal{S}$ 
is a cellular algebra in the sense of \cite{GL}. Given $\lambdal \in \Pilm$, one defines as in \cite[Definition 6.13]{DJM} a right 
$\mathcal{S}$-module $W(\lambdal)$ which is a free $R$-module of finite rank, called \emph{Weyl module}. Since $\mathcal{S}$ is cellular, 
$W(\lambdal)$ is naturally equipped with a symmetric bilinear form $\langle \cdot ,\cdot \rangle$. Set 
\begin{equation}
L(\lambdal):=W(\lambdal) / \mbox{rad} \,  W(\lambdal),
\end{equation}
where $\mbox{rad} \,  W(\lambdal)$ is the radical of the bilinear form $\langle \cdot ,\cdot \rangle$. Assume temporarily that $R$ 
is a field. By \cite[Corollary 6.18]{DJM}, $\mathcal{S}$ is a quasi-hereditary algebra, so the theory of cellular algebras of \cite{GL} shows that 
$\set{L(\lambdal) \mid \lambdal \in \Pilm}$ is a complete set of non-isomorphic irreducible $\mathcal{S}$-modules (see \cite[Theorem 6.16]{DJM}).
This implies that $\mathcal{R}_0(\mathcal{S})$, the Grothendieck group of finitely-generated $\mathcal{S}$-modules, is a free $\mZ$-module with 
basis $\set{[L(\lambdal)] \mid \lambdal \in \Pilm}$. \\ 

\! \,From now on, we assume that $R$ is a local ring, with unique maximal ideal $\wp$. Let $\nu_{\wp}$ be the corresponding $\wp$-adic valuation 
map. Let $K$ be the field of fractions of $R$ and extend $\nu_{\wp}$ to $K$ in the natural way. Let $F=R/\wp R$ be the residue field, so $(R,K,F)$ 
is a modular system. If $M$ is a right $R$-module, we denote by $M_F = M \otimes_R F$ the specialized module and denote similarly by
$M_K = M \otimes_R K$ the corresponding module defined over 
$K$. We shall use this notation for Weyl modules and for $\mathcal{S}$ itself. \\ 

\begin{definition}[\cite{Jan}, see also \cite{AM}]
Let $M$ be an $R$-module equipped with a symmetric bilinear form  $\langle \cdot ,\cdot \rangle$. For all $i \in \mN$, set 
\begin{equation}
M(i):=\set{u \in M \ | \ \forall \ v \in M, \ \nu_{\wp}(\langle u,v \rangle) \geq i}.
\end{equation}
The \emph{Jantzen filtration} of $M$ is the sequence
\begin{equation}
M_F = M_F(0) \supset M_F(1) \supset \cdots,
\end{equation} 
where $M_F(i):=(M(i)+\wp M)/\wp M$. \fini
\end{definition}

\noindent Note that in the definition above, we have in particular $M_F(1) = \mbox{rad} \, M_F$. Moreover, if $M$ is free of finite 
rank (as an $R$-module), then we have $M_F(i) = \set{0}$ for $i$ large enough. \\

The following theorem was proved by James and Mathas (see \cite[Theorem 4.3]{JM}).

\begin{thm}[the Jantzen sum formula] \label{thm_Jantzen_formula} \ \\ 
Assume that $\mathcal{S}_K$ is semisimple. Then in the Grothendieck group 
$\mathcal{R}_0(\mathcal{S}_F)$, we have for all $\lambdal \in \Pilm$: 
\begin{equation} 
\ds \sum_{i>0}{[W_F(\lambdal;i)]} = \sum_{\mul \in \Pilm}{\nu_{\wp}(g_{\lambdal,\mul})\, [W_F(\mul)]}.
\end{equation}
Here, the $g_{\lambdal,\mul} \in R$ are factors of some Gram determinants (see \cite[Definitions 3.1, 3.36 and Corollary 3.38]{JM}). \cqfd
\end{thm}

\begin{remark} The condition of semisimplicity of $\mathcal{S}_K$ is stated in \cite[Theorem 4.3]{JM} in terms of the Poincar\'e polynomial for 
$\mathcal{H}_R$, which is defined in \cite[Definition 3.40]{JM}. \fini
\end{remark}

James and Mathas \cite{JM} showed that only multi-partitions $\mul \in \Pilm$ such that $\mul \lhd \lambdal$ contribute to the right hand-side of
Theorem \ref{thm_Jantzen_formula} (the definition of the dominance ordering $\lhd$ is recalled in Definition \ref{def_dominance_ordering}). 
They have given a combinatorial expression of $\nu_{\wp}(g_{\lambdal,\mul})$ in terms of ribbons 
contained in diagrams of $l$-multi-partitions. However, this combinatorial expression makes sense even if
$\lambdal$ does not dominate $\mul$. We will therefore introduce in Section \ref{def_J} a matrix 
$J=\bigl(j_{\lambdal,\mul} \bigr)_{\lambdal, \mul \in \Pilm}$ whose entries are these combinatorial expressions \underline{without restriction on 
the pair $(\lambdal, \mul)$}. More precisely, our indexing is chosen so that
\begin{equation} \label{eq_matrix_J_Jantzen_sum_formula}
j_{\lambdal^{\dagger},\mul^{\dagger}}=\nu_{\wp}(g_{\lambdal,\mul}) \qquad \mbox{if} \quad \mul \lhd \lambdal,
\end{equation}
where the sign $\dagger$ denotes the conjugation of multi-partitions (see (\ref{equation_conjugate})). We are forced to use conjugates here because 
the indexation from \cite{JM} for the rows and columns of decomposition matrices is not compatible with the indexation from \cite{U2} for the rows 
and columns of transition matrices for Uglov's canonical bases.

Now, let $\leqslant$ be an arbitrary partial ordering on $\Pilm$ and write $\lambdal < \mul$ if 
$\lambdal \leqslant \mul$ and $\lambdal \neq \mul$ ($\lambdal,\, \mul \in \Pilm$). 
Define a matrix $J^{<}=\bigl(j^{<}_{\lambdal,\mul} \bigr)_{\lambdal,\, \mul \in \Pilm}$ by the formula
\begin{equation} \label{eq_def_J_inf}
j^{<}_{\lambdal,\mul}:= \left\{ 
\begin{array}{cc}
j_{\lambdal,\mul} & \mbox{if } \lambdal < \mul \\[1mm]
0 & \mbox{otherwise}
\end{array}
\right. \qquad (\lambdal,\, \mul \in \Pilm) .
\end{equation}
If we take $\leqslant \ = \ \unlhd$, then we get a matrix $J^{\lhd}$ whose entries are,
up to conjugation of multi-partitions, the $\nu_{\wp}(g_{\lambdal,\mul})$'s of \cite{JM}. \\ 
   
We now derive a matrix identity equivalent to the Jantzen sum formula. 

\begin{definition} \label{definition_Dq}
Let $D(q)=\bigl(d_{\lambdal,\mul}(q) \bigr)_ {\lambdal, \mul \in \Pilm}$ be the matrix defined by
\begin{equation}
d_{\lambdal,\mul}(q) := \sum_{i \geq 0}{\big[W_F(\lambdal^{\dagger};i)/W_F(\lambdal^{\dagger};i+1) : L_F(\mul^{\dagger})\big]  q^i} 
\in \mN[q] \qquad (\lambdal, \mul \in \Pilm).
\end{equation} 
\fini
\end{definition}

\noindent Note that $d_{\lambdal^{\dagger},\mul^{\dagger}}(1)$ is 
equal to the multiplicity of $L_F(\mul)$ as a composition factor of $W_F(\lambdal)$, so up to conjugation of  
multi-partitions (which amounts to reindexing the rows and columns of the matrix), $D(1)$ is the usual  
decomposition matrix of $\mathcal{S}_F$.

\begin{lemma} \label{lemma_Jantzen_formula}
Let $M=\bigl(m_{\lambdal,\mul} \bigr)_{\lambdal,\mul \in \Pilm}$ be a matrix with integer entries. Then the following statements are 
equivalent:
 
\begin{itemize}
\item[\rm (i)] In $\mathcal{R}_0(\mathcal{S}_F)$, we have for all $\lambdal \in \Pilm\,$: 
$\ds \sum_{i>0}{[W_F(\lambdal^{\dagger};i)]} = \sum_{\nul \in \Pilm}{m_{\lambdal,\nul}[W_F(\nul^{\dagger})]},$
\item[\rm (ii)] $D'(1)=MD(1)$.
\end{itemize}
\end{lemma}

\begin{demo} 
Let $\lambdal \in \Pilm$. Since $\set{\big[L_F(\mul^{\dagger})\big] \ \big| \ \mul \in \Pilm}$ is a $\mZ$-basis of $\mathcal{R}_0(\mathcal{S}_F)$, 
we have on the one hand: 
\bigskip
$$\begin{array}{rcl} \ds \sum_{i>0}{\big[W_F(\lambdal^{\dagger}\,;\,i)\big]} &=&
\ds \sum_{i>0}{\sum_{\mul \in \Pilm}{\big[W_F(\lambdal^{\dagger}\,;\,i) : L_F(\mul^{\dagger})\big] 
\big[L_F(\mul^{\dagger})\big]}} \\[5mm]
&=& \ds \sum_{\mul \in \Pilm}{\Bigl({ \sum_{i>0}{\sum_{j \geq i}
{\big[W_F(\lambdal^{\dagger}\,;\,j)\,/\,W_F(\lambdal^{\dagger}\,;\,j+1) : L_F(\mul^{\dagger})\big]}}}\Bigr)
\big[L_F(\mul^{\dagger})\big]} \\[5mm]
&=& \ds \sum_{\mul \in \Pilm}
{\Bigl({ \sum_{j>0}{\sum_{0 < i \leq j}{\big[W_F(\lambdal^{\dagger}\,;\,j)\,/\,W_F(\lambdal^{\dagger}\,;\,j+1) : 
L_F(\mul^{\dagger})\big]}}}\Bigr)\big[L_F(\mul^{\dagger})\big]} \\[5mm]
&=& \ds \sum_{\mul \in \Pilm}{d'_{\lambdal,\mul}(1)\big[L_F(\mul^{\dagger})\big]}.
\end{array}$$

On the other hand, we have 
\bigskip
$$\begin{array}{rcl} \ds \sum_{\nul \in \Pilm}{m_{\lambdal,\nul}\big[W_F(\nul^{\dagger})\big]} &=&
 \ds \sum_{\mul, \, \nul \in \Pilm}{m_{\lambdal,\nul}}
 \big[W_F(\nul^{\dagger}):L_F(\mul^{\dagger})\big]\big[L_F(\mul^{\dagger})\big] \\[5mm]
&=& \ds \sum_{\mul \in \Pilm}
{\Bigl({\sum_{\nul \in \Pilm}{m_{\lambdal,\nul}\big[W_F(\nul^{\dagger}):L_F(\mul^{\dagger})\big]}}\Bigr) 
\big[L_F(\mul^{\dagger})\big]} \\[5mm] 
&=& \ds \sum_{\mul \in \Pilm}
{\Bigl({\sum_{\nul \in \Pilm}{m_{\lambdal,\nul}d_{\nul,\mul}(1)}}\Bigr) \big[L_F(\mul^{\dagger})\big]}.
\end{array}$$ 
since the $\big[L_F(\mul^{\dagger})\big]$, $\mul \in \Pilm$ are linearly independent, the lemma follows.
\end{demo}

The Jantzen sum formula as stated in Theorem \ref{thm_Jantzen_formula}, together with (\ref{eq_matrix_J_Jantzen_sum_formula}) and Lemma 
\ref{lemma_Jantzen_formula}, implies the following result.

\begin{cor} \label{cor_Jantzen_formula} Assume that $\mathcal{S}_K$ is semisimple. Then with the notation above,
we have \begin{equation} D'(1)=J^{\lhd}D(1). \end{equation} \cqfd
\end{cor} 

\subsection{Statement of theorems} 

In this section, we state an important conjecture for computing the decomposition matrix of the cyclotomic $v$-Schur algebra
defined over $\mC$, with parameters equal to arbitrary powers of a primitive $n$-th root of unity. This conjecture is supported by Theorem 
\ref{thm_Delta_J}.  

\subsubsection{Choice of parameters} \label{section_choice_parameters}

\! \,Fix $(\mathtt{r}_1,\ldots,\mathtt{r}_l) \in (\mZ/n\mZ)^l$. We shall define a modular system $(R,K,F)$ with parameters such that 
the specialized cyclotomic $v$-Schur algebra $\mathcal{S}_F$ is
$\mathcal{S}_{\mC,m}(\zeta;\zeta^{\mathtt{r}_1},\ldots,\zeta^{\mathtt{r}_l})$ with $\zeta:=\exp(\frac{2i\pi}{n})$. \\

We first define a modular system $(R,K,F)$ as follows. Let $\widehat{R}=\mC[x,\,x^{-1}]$ be the ring of Laurent polynomials in one indeterminate
over the field $\mC$. Let
\begin{equation}
\begin{array}{c}
\xi:=\exp \Bigl(\frac{2i\pi}{nl} \Bigr) \in \mC, \quad \wp:=(x-\xi), \quad R:=\mC[x,\,x^{-1}]_{\wp}, \\[3mm]
 K:=\mC(x) \quad \mbox{and} \quad F:=R/\wp R \simeq \mC,
\end{array}
\end{equation}
that is, $\wp$ is the prime ideal in $\widehat{R}$ spanned by $x-\xi$ with $\xi$ a primitive complex $nl$-th root of unity, $R$ is
the localized ring of $\widehat{R}$ at $\wp$, $K$ is the field of fractions of $R$ and $F$ is the residue field. \\

\! \,Following \cite{U2}, we fix an $l$-tuple $\esselle$ in
\begin{equation}
\mathcal{L}(\mathtt{r}_1,\ldots,\mathtt{r}_l):=\set{(s_1,\ldots,s_l) \in \mZ^l \mid \forall \, 1 \leq d \leq l,\, \mathtt{r}_d = s_d \bmod n}.    
\end{equation}

\noindent Such an $l$-tuple is called a \emph{multi-charge}. The multi-charge $\esselle$ parametrizes a so-called \emph{($q$-deformed) Fock space of 
level $l$}, denoted by $\Fq$ (see Section \ref{section_fock_space}). Note that for a given 
$(\mathtt{r}_1,\ldots,\mathtt{r}_l) \in (\mZ/n\mZ)^l$ we have an infinite choice of Fock spaces $\Fq$ such that 
$\esselle$ is in $\mathcal{L}(\mathtt{r}_1,\ldots,\mathtt{r}_l)$. \\

We now describe the choice of parameters for the cyclotomic $v$-Schur algebra $\mathcal{S}$. These parameters are similar to those used in 
\cite{Jac} for Ariki-Koike algebras. They depend on $n$, $l$ and on the multi-charge 
$(s_1,\ldots,s_l) \in \mathcal{L}(\mathtt{r}_1,\ldots,\mathtt{r}_l)$ that we have fixed. Put
\begin{equation}
v:=x^l \qquad  \mbox{and} \qquad u_d:= \xi^{nd} x^{ls_d-nd} \qquad (1 \leq d \leq l).
\end{equation}

\noindent Note that we have $\mathcal{S}_{F}=\mathcal{S}_{\mC,m}(\zeta;\zeta^{\mathtt{r}_1},\ldots,\zeta^{\mathtt{r}_l})$ with 
$\zeta:=\exp(\frac{2i\pi}{n})$. Note also that the algebra $\mathcal{S}_{K,m}(v;u_1,\ldots,u_l)$ is semisimple. Indeed, specializing $x$ at $1$ sends 
$\mathcal{H}_{K,m}(v;u_1,\ldots,u_l)$ on the semisimple group algebra $\mC G(l,1,m)$, so by the Tits deformation argument \cite{A1}, the algebra 
$\mathcal{H}_{K,m}(v;u_1,\ldots,u_l)$ is semisimple and so is $\mathcal{S}_{K,m}(v;u_1,\ldots,u_l)$. Therefore, the Jantzen sum formula
(see Theorem \ref{thm_Jantzen_formula}) applies in our case. This leads in particular to the definition of a matrix $J^{\lhd}$ 
(see Section \ref{def_J}). 

\subsubsection{Main result} \label{section_our_results}

\! \,Following \cite{U2}, let $\esselle \in \mathcal{L}(\mathtt{r}_1,\ldots,\mathtt{r}_l)$ and $\Fq$ be the corresponding Fock space of level $l$ 
(see Section \ref{section_fock_space}). As a vector space, $\Fq$ has a natural basis $\set{ |\lambdal,\esselle \rangle \mid \lambdal \in \Pi^l}$ 
and a canonical basis $\set{\mathcal{G}(\lambdal,\esselle) \mid \lambdal \in \Pi^l}$ indexed by $l$-multi-partitions. Let $\Fq_{m}$ be the subspace 
of $\Fq$ spanned by the $|\lambdal,\esselle \rangle$'s, $\lambdal \in \Pilm$. Let $A(q)$ be the matrix of the involution $\barre$ of $\Fq_{m}$ with 
respect to the standard basis, and let $\Delta(q)$ be the transition matrix between the standard basis and the canonical basis of $\Fq_{m}$ 
(see Sections \ref{section_fock_involution} and \ref{section_canonical_basis}). Still following \cite{U2}, we associate to $\esselle$ an ordering
$\prec$ (see Definition \ref{def_prec}). By (\ref{eq_def_J_inf}) we get a matrix $J^{\prec}$. 

\begin{thm} \label{thm_Delta_J} Let $\esselle \in \mathcal{L}(\mathtt{r}_1,\ldots,\mathtt{r}_l)$. Then with the notation above, we have
\begin{equation} \Delta'(1)=J^{\prec} \Delta(1). 
\end{equation} \cqfd
\end{thm}

\begin{example} \label{example_1_thm_Delta_J}
Take $n=3$, $l=2$, $\esselle=(1,0)$ and $m=3$. Then we have on the one hand \\

$$J^{\prec}= \begin{array}{ll}
 \left( \begin{array}{cccccccccc}
 0 & . & . & . & . & . & . & . & . & .  \\
 0 & 0 & . & . & . & . & . & . & . & .  \\
 0 & 0 & 0 & . & . & . & . & . & . & .  \\ 
 0 & 0 & 0 & 0 & . & . & . & . & . & .  \\
 0 & 0 & 1 & 0 & 0 & . & . & . & . & .  \\
 0 & 1 & 1 & 0 & 0 & 0 & . & . & . & .  \\
 0 &-1 & 0 & 0 & 1 & 1 & 0 & . & . & .  \\
 0 & 1 &-1 & 0 & 1 &-1 & 1 & 0 & . & .  \\
 0 & 1 &-1 & 0 & 0 & 1 & 0 & 0 & 0 & .  \\
 0 &-1 & 0 & 0 & -1& 0 & 1 & 0 & 1 & 0  \\
 \end{array} \right)
&
 \begin{array}{l} 
 \bigl( (1,1),(1) \bigr) \\
 \bigl( (3),\emptyset \bigr) \\
 \bigl( \emptyset,(3) \bigr) \\
 \bigl( (1),(2) \bigr) \\
 \bigl( \emptyset,(2,1) \bigr) \\
 \bigl( (2),(1) \bigr) \\
 \bigl( (1),(1,1) \bigr) \\
 \bigl( \emptyset,(1,1,1) \bigr) \\
 \bigl( (2,1),\emptyset \bigr) \\
 \bigl( (1,1,1),\emptyset \bigr) \\
 \end{array}
\end{array} ,$$

\vspace{5mm}
\noindent where dots over the main diagonal stand for zero entries. The $l$-multi-partitions of $m$ which index the bases of $\Fq_m$ 
are ordered decreasingly with respect to a total ordering finer than $\prec$ and they are displayed in the column located on the right 
of the matrix $J^{\prec}$. On the other hand, we compute $\Delta(q)$ using Uglov's algorithm (see \cite{U2}). If we keep the same 
ordering for the rows and the columns of $\Delta(q)$, we get the following matrix. 

$$\Delta(q)= \begin{array}{ll} 
 \left( \begin{array}{cccccccccc}
 1 & . & . & . & . & . & . & . & . & . \\
0 &1 & . & . & . & . & . & . & . & . \\
0 &0 & 1 & . & . & . & . & . & . & . \\
0 &0 & 0 & 1 & . & . & . & . & . & . \\
0 &0 & q & 0 & 1 & . & . & . & . & . \\
0 &q & q & 0 & 0 & 1 & . & . & . & . \\
0 &0 & q^2 & 0 & q & q & 1 & . & . & . \\
0 &0 & 0 & 0 & q^2 & 0 & q & 1 & . & . \\ 
0 &q^2 & 0 & 0 & 0 & q & 0 & 0 & 1 & . \\
0 &0 & 0 & 0 & 0 & q^2 & q & 0 & q & 1 \\
\end{array} \right)
&
 \begin{array}{l} 
 \bigl( (1,1),(1) \bigr) \\
 \bigl( (3),\emptyset \bigr) \\
 \bigl( \emptyset,(3) \bigr) \\
 \bigl( (1),(2) \bigr) \\
 \bigl( \emptyset,(2,1) \bigr) \\
 \bigl( (2),(1) \bigr) \\
 \bigl( (1),(1,1) \bigr) \\
 \bigl( \emptyset,(1,1,1) \bigr) \\
 \bigl( (2,1),\emptyset \bigr) \\
 \bigl( (1,1,1),\emptyset \bigr) \\
 \end{array}
\end{array}.$$

\vspace{5mm}
\noindent It is easy to check that $\Delta'(1)=J^{\prec} \Delta(1)$. \fini 
\end{example} 

\begin{example} \label{example_2_thm_Delta_J}
Take $n=3$, $l=2$, $\esselle=(4,-3)$ and $m=3$. Write the rows and the columns of the following matrices with respect to a total ordering finer 
than $\prec$. Then \\

$$J^{\prec}= \begin{array}{ll}
 \left( \begin{array}{cccccccccc}
 0 & . & . & . & . & . & . & . & . & .  \\
 0 & 0 & . & . & . & . & . & . & . & .  \\
 0 & 1 & 0 & . & . & . & . & . & . & .  \\ 
 0 & 1 & 1 & 0 & . & . & . & . & . & .  \\
 0 & -1& 1 & 0 & 0 & . & . & . & . & .  \\
 0 & 0 & 0 & 0 & 0 & 0 & . & . & . & .  \\
 0 &-1 & 0 & 1 & 1 & 0 & 0 & . & . & .  \\
 0 & 0 &-1 & 1 & 0 & 0 & 0 & 0 & . & .  \\
 0 & 0 & 0 & 0 & -1& 0 & 1 & 1 & 0 & .  \\
 0 & 1 & 0 & -1& 0 & 0 & 1 & -1& 1 & 0  \\
 \end{array} \right)
&
 \begin{array}{l} 
 \bigl( (1,1),(1) \bigr) \\
 \bigl( (3),\emptyset \bigr) \\
 \bigl( (2,1),\emptyset \bigr) \\
 \bigl( (2),(1) \bigr) \\
 \bigl( (1,1,1),\emptyset \bigr) \\
 \bigl( (1),(2) \bigr) \\
 \bigl( (1),(1,1) \bigr) \\
 \bigl( \emptyset,(3) \bigr) \\
 \bigl( \emptyset,(2,1) \bigr) \\ 
 \bigl( \emptyset,(1,1,1) \bigr) \\
 \end{array}
\end{array}$$

\vspace{5mm}
and \\

$$\Delta(q) = \begin{array}{ll}
 \left( \begin{array}{cccccccccc}
 1  & .  & .   & .   & .  & .   & .   & .  & . & .  \\
 0 	& 1  & .   & .   & .  & .   & .   & .  & . & .  \\
 0 	& q	 & 1   & .   & .  & .   & .   & .  & . & .  \\ 
 0 & q^2 & q   & 1   &  . & .   & .   & .  & . & .  \\
 0 & 0	 & q   & 0   & 1  & .   & .   & .  & . & .  \\
 0 & 0	 & 0   & 0   & 0  &  1  & .   & .  & . & .  \\
 0 & 0   & q^2 & q   & q  &  0  & 1   & .  & . & .  \\
 0 & 0   & 0   & q   & 0  &  0  & 0   & 1  & . & .  \\
 0 & 0   & 0   & q^2 & 0  &  0  & q   & q  & 1 & .  \\
 0 & 0   & 0   & 0   & q  & 0   & q^2 & 0  & q & 1  \\
 \end{array} \right)
&
 \begin{array}{l} 
 \bigl( (1,1),(1) \bigr) \\
 \bigl( (3),\emptyset \bigr) \\
 \bigl( (2,1),\emptyset \bigr) \\
 \bigl( (2),(1) \bigr) \\
 \bigl( (1,1,1),\emptyset \bigr) \\
 \bigl( (1),(2) \bigr) \\
 \bigl( (1),(1,1) \bigr) \\
 \bigl( \emptyset,(3) \bigr) \\
 \bigl( \emptyset,(2,1) \bigr) \\ 
 \bigl( \emptyset,(1,1,1) \bigr) \\
 \end{array}
\end{array}.$$ 

\vspace{5mm}
\noindent Again, one can check that $\Delta'(1)=J^{\prec} \Delta(1)$. \fini 
\end{example}

\bigskip

Theorem \ref{thm_Delta_J} is equivalent to the following: \\

\begin{thm} \label{thm_A_J} With the notation of Theorem \ref{thm_Delta_J}, we have
\begin{equation} A'(1)=2J^{\prec}.
\end{equation}
\end{thm}

\noindent \emph{Proof of the equivalence of Theorems \ref{thm_Delta_J} and \ref{thm_A_J}.} Since the canonical basis is invariant under the 
$\barre$ involution, we have $\Delta(q)=A(q)\Delta(q^{-1})$. Taking derivatives at $q=1$ yields $\Delta'(1)=A'(1)\Delta(1)-A(1)\Delta'(1)$.
Since $A(1)$ is the identity matrix, we get $2\Delta'(1)=A'(1)\Delta(1)$. As a consequence, Theorem \ref{thm_A_J} implies Theorem 
\ref{thm_Delta_J}. Since $\Delta(1)$ is unitriangular, hence invertible, the converse follows. \cqfd \\ 

We prove Theorem \ref{thm_A_J} in Part C. Our proof is similar to the proof of \cite{Ry} in the level one case.
However the higher-level case is significantly more complicated and involves the discussion of many cases (see Section \ref{proof_thm_A_J}). 

\subsubsection{A conjecture for the decomposition matrix of $\mathcal{S}$}

Choose the parameters as in Section \ref{section_choice_parameters}. Guided by the formal analogy between Theorem \ref{thm_Delta_J} on one hand,
and the rephrasing of the Jantzen sum formula given in Corollary \ref{cor_Jantzen_formula} on the other hand, we may wonder if
for some $\esselle \in \mathcal{L}(\mathtt{r}_1,\ldots,\mathtt{r}_l)$, the corresponding matrix $J^{\prec}$ coincides with the matrix $J^{\lhd}$ 
coming from the Jantzen sum formula. This leads to the following definition and conjecture.

\begin{definition} \label{def_dominant}
Let $M \in \mN$. We say that $\esselle \in \mathcal{L}(\mathtt{r}_1,\ldots,\mathtt{r}_l)$ is \emph{$M$-dominant} if for all 
$1 \leq d \leq l-1$, we have 
\begin{equation}
s_{d+1}-s_d \geq M.
\end{equation} \fini
\end{definition}

\noindent The point is that if $\esselle$ is $m$-dominant, then we have $J^{\prec}=J^{\lhd}$ (see Proposition \ref{expr_J_dominant}). 

\begin{conj} \label{main_conj} Assume that $\esselle \in \mathcal{L}(\mathtt{r}_1,\ldots,\mathtt{r}_l)$ is $m$-dominant. Let $D(q)$ be the 
$q$-analogue of the decomposition matrix of $\mathcal{S}$ defined in Definition \ref{definition_Dq} with our choice of parameters given in Section 
\ref{section_choice_parameters}. Then we have 
\begin{equation}
D(q)=\Delta(q).
\end{equation} \cqfd
\end{conj}

\noindent\,\!If we put $q=1$ in Conjecture \ref{main_conj}, we thus get an algorithm for computing the decomposition matrix of
$\mathcal{S}_{\mC,m}(\zeta;\zeta^{\mathtt{r}_1},\ldots,\zeta^{\mathtt{r}_l})$ with $\zeta:=\exp(\frac{2i\pi}{n})$.

\begin{remark} The assumption of $m$-dominance is necessary in Conjecture \ref{main_conj}. Indeed, while the decomposition matrix $D(1)$ only
depends on the sequence $(\mathtt{r}_1,\ldots,\mathtt{r}_l)$ of the residues modulo $n$ of the multi-charge $\esselle$,
the matrix $\Delta(1)$ actually depends on $\esselle$ itself. For example, take $n=3$, $l=2$ and $m=3$. Then the multi-charges $(1,0)$ and $(4,-3)$ 
are both in $\mathcal{L}(1,0)$, but the corresponding matrices $\Delta(1)$ do not have the same number of zero entries
(see Examples \ref{example_1_thm_Delta_J} and \ref{example_2_thm_Delta_J}). \fini \\ 
\end{remark}   

\begin{remark} Conjecture \ref{main_conj} suggests that the matrix $\Delta(1)$ should not depend of the choice of the multi-charge 
$\esselle \in \mathcal{L}(\mathtt{r}_1,\ldots,\mathtt{r}_l)$ provided it is $M$-dominant for $M$ large enough. This statement is proved in
\cite[Th\'eorème 4.30]{Y}), where an explicit value of $M$ is given. However, the fact that we might take $M=m$ here is still conjectural. \fini
\end{remark} 

\begin{example} Set $n=3$, $l=2$, $(\mathtt{r}_1,\mathtt{r}_2)=(1,0)$ and $m=3$. Then the specialized cyclotomic $v$-Schur algebra is 
$\mathcal{S}_{\mC,3} \Bigl(e^{\frac{2i \pi}{3}};e^{\frac{2i \pi}{3}},1 \Bigr)$. Take $\esselle=(4,-3)$, so 
$\esselle \in \mathcal{L}(\mathtt{r}_1,\ldots,\mathtt{r}_l)$
is $m$-dominant. According to Conjecture \ref{main_conj}, we expect $D(q)$ be equal to \\  

$$\Delta(q) = \begin{array}{ll}
 \left( \begin{array}{cccccccccc}
 1  & .  & .   & .   & .  & .   & .   & .  & . & .  \\
 0 	& 1  & .   & .   & .  & .   & .   & .  & . & .  \\
 0 	& q	 & 1   & .   & .  & .   & .   & .  & . & .  \\ 
 0 & q^2 & q   & 1   &  . & .   & .   & .  & . & .  \\
 0 & 0	 & q   & 0   & 1  & .   & .   & .  & . & .  \\
 0 & 0	 & 0   & 0   & 0  &  1  & .   & .  & . & .  \\
 0 & 0   & q^2 & q   & q  &  0  & 1   & .  & . & .  \\
 0 & 0   & 0   & q   & 0  &  0  & 0   & 1  & . & .  \\
 0 & 0   & 0   & q^2 & 0  &  0  & q   & q  & 1 & .  \\
 0 & 0   & 0   & 0   & q  & 0   & q^2 & 0  & q & 1  \\
 \end{array} \right)
&
 \begin{array}{l} 
 \bigl( (1,1),(1) \bigr) \\
 \bigl( (3),\emptyset \bigr) \\
 \bigl( (2,1),\emptyset \bigr) \\
 \bigl( (2),(1) \bigr) \\
 \bigl( (1,1,1),\emptyset \bigr) \\
 \bigl( (1),(2) \bigr) \\
 \bigl( (1),(1,1) \bigr) \\
 \bigl( \emptyset,(3) \bigr) \\
 \bigl( \emptyset,(2,1) \bigr) \\ 
 \bigl( \emptyset,(1,1,1) \bigr) \\
 \end{array}
\end{array}$$ (see Example \ref{example_2_thm_Delta_J}). \fini
\end{example}

If we no longer assume that $\esselle$ is $m$-dominant, then we expect $\Delta(q)$ be equal to a $q$-analogue of the decomposition matrix of a 
quasi-hereditary covering (in the sense of Rouquier, see \cite{Ro}) of the Ariki-Koike algebra $\mathcal{H}$. This covering, depending on 
$\esselle$, could come from a rational Cherednik algebra through the Knizhnik-Zamolodchikov functor \cite{GGOR}. It should be Morita-equivalent 
to the cyclotomic $v$-Schur algebra of \cite{DJM} if $\esselle$ is $m$-dominant.   

\vspace{5mm}
\begin{center}
\scshape \Large PART B: Tools for the proof of Theorem \ref{thm_A_J}
\end{center}
\vspace{5mm}

The next two sections recall some results about combinatorics of partitions and multi-partitions on the one hand and higher-level Fock spaces on
the other hand; all of them will be used in the proof of Theorem \ref{thm_A_J}. However, there are no new results here, so the reader 
familiar with these two topics may skip this part and come back to it later in order to get the needed definitions and notation.   

\section{Combinatorics of partitions and multi-partitions}
\subsection{Definitions} 

We give here all the basic definitions about partitions and multi-partitions that we need later; our main reference is $\cite{Mac}$. Let 
$r \in \mN$. A \emph{partition of $r$} is a sequence of integers $\lambda=(\lambda_1,\lambda_2,\ldots,\lambda_N)$ such that
$\lambda_1 \geq \lambda_2 \geq \ldots \geq \lambda_N \geq 0$ and $\lambda_1 + \cdots + \lambda_N = r$. Each nonzero $\lambda_i$ is called a
\emph{part} of $\lambda$. The sum of all the parts of $\lambda$ is denoted by $|\lambda|$. We identify two partitions differing only by a 
tail of zeroes and write sometimes partitions as sequences of integers with an infinite tail of zeroes. The only partition of $0$ is denoted by 
$\emptyset$. The \emph{conjugate} of the partition $\lambda$ is the partition $\lambda^{\dagger}$ defined by 
\begin{equation} \label{equation_conjugate}
\lambda^{\dagger}_i:=\sharp \set{j \ | \ \lambda_j \geq i} \qquad (i \geq 1) ;
\end{equation}
for example, the conjugate of $(4,3,3,2,1)$ is $(5,4,3,1)$. \\

An \emph{$N$-multi-partition of $r$} is an $N$-tuple of partitions of integers summing up to $r$. Let 
$\boldsymbol{\lambda} = (\lambda^{(1)},\ldots,\lambda^{(N)})$ be an $N$-multi-partition. The \emph{conjugate} of $\boldsymbol{\lambda}$ is the 
multi-partition $\boldsymbol{\lambda}^{\dagger}:=\big((\lambda^{(N)})^{\dagger},\ldots,(\lambda^{(1)})^{\dagger} \big)$. 
For $1 \leq b \leq N$, write $\lambda^{(b)} = (\lambda^{(b)}_1,  \lambda^{(b)}_2,  \ldots)$ the parts of $\lambda^{(b)}$. 
The \emph{Young diagram} of $\boldsymbol{\lambda}$ is the set 
\begin{equation}
\set{(i,j,b) \in \mN^* \times \mN^* \times \interv{1}{N} \ | \ 1 \leq j \ \leq \lambda^{(b)}_i} ,
\end{equation}
whose elements are called \emph{boxes} or \emph{nodes} of $\boldsymbol{\lambda}$. If $N=1$, namely, if $\boldsymbol{\lambda}$ is a partition, 
we drop the third component in the symbol $(i,j,b)$ of a node of $\boldsymbol{\lambda}$. From now on we identify an $N$-multi-partition with
its Young diagram. We extend the notation $|\lambda|$ in a natural way for multi-partitions and define the \emph{dominance ordering} on 
multi-partitions as follows. 
\begin{definition} \label{def_dominance_ordering}
Let $\boldsymbol{\lambda}$ and $\boldsymbol{\mu}$ be two $N$-multi-partitions. We say that
$\boldsymbol{\mu}$ \emph{dominates} $\boldsymbol{\lambda}$ and write $\boldsymbol{\lambda} \unlhd \boldsymbol{\mu}$ if
\begin{equation} \label{condition1_dominance_ordering}
|\boldsymbol{\lambda}|=|\boldsymbol{\mu}|
\end{equation}
 and for all $k \geq 0$, $1 \leq b \leq N$, we have 
\begin{equation}
\ds \sum_{i=1}^{b-1} |\lambda^{(i)}| + \sum_{j=1}^{k} \lambda^{(b)}_j \leq \ds \sum_{i=1}^{b-1} |\mu^{(i)}| + \sum_{j=1}^{k} \mu^{(b)}_j.
\end{equation}
\noindent Write $\boldsymbol{\lambda} \lhd \boldsymbol{\mu}$ if $\boldsymbol{\lambda} \unlhd \boldsymbol{\mu}$ and 
$\boldsymbol{\lambda} \neq \boldsymbol{\mu}$. \fini
\end{definition}
If $\lambda, \mu \in \Pi$ are two partitions, write $\lambda \subset \mu$ if the diagram of $\lambda$ is contained in the diagram of $\mu$,
and the set-theoretic difference is called a $\emph{skew diagram}$; we denote it by $\mu/\lambda$. A \emph{path} in the skew diagram $\theta$ is
a sequence of boxes $(\gamma_1,\ldots,\gamma_N) \in \theta^N$ such that for all $1 \leq i \leq N-1$, $\gamma_i$ and $\gamma_{i+1}$ have one common
side. We say that $\theta$ is \emph{connected} if given any two boxes $\gamma, \gamma' \in \theta$, there exists a path within $\theta$ connecting
$\gamma$ to $\gamma'$. A \emph{ribbon} is a connected skew diagram that contains no $2 \times 2$ block of boxes. Let $\rho$ be a ribbon. The   
\emph{head} (\resp \emph{tail}) of $\rho$ is the node $\gamma=(i,j) \in \rho$ such that $j-i$ is minimal (\resp maximal); we denote this node by
$\hd(\rho)$ (\resp $\tl(\rho)$). If $\hd(\rho)=(i,j)$ and $\tl(\rho)=(i',j')$, the \emph{height} of $\rho$ is the integer 
$\haut(\rho):=i-i' \in \mN$. Finally, the \emph{length} of $\rho$ is the number of boxes contained in $\rho$; we denote it by $\ell(\rho)$.   

\begin{example} On Figure \ref{fig:2} (see Section \ref{def_J}), the set of white squares represents the partition $(4,1)$ ; $\rho$, 
$\rho'$ and $\rho''$ are three ribbons of respective heights $2$, $1$ and $0$ and of respective lengths $4$, $4$ and $3$. \fini
\end{example}

A \emph{charged $N$-multi-partition} is an element of $\Pi^N \times \mZ^N$. If $(\boldsymbol{\lambda},\textbf{s}) \in \Pi^N \times \mZ^N$ is a
charged multi-partition and $\textbf{s}=(s_1,\ldots,s_N)$, the \emph{content} of the node $\gamma=(i,j,b) \in \boldsymbol{\lambda}$ is the 
integer 
\begin{equation}
\cont(\gamma) := s_b+j-i.
\end{equation}
If $M \in \mN^{*}$, the \emph{residue modulo} $M$ of $\gamma$ is  
\begin{equation}
\res_M(\gamma):=\cont(\gamma) \bmod M \in \mZ/M\mZ.
\end{equation} 
\! \,For all $i \in \mZ$, set
\begin{equation} \label{notation_nb_i_nodes}
N_i(\boldsymbol{\lambda}):= \sharp \set{\gamma \in \boldsymbol{\lambda} \ | \ \res_n(\gamma)=i \bmod n}\, ;
\end{equation}
this number depends on the multi-charge $\textbf{s}$. Define in a similar way $N_i(\theta)$ if $\theta$ is a skew diagram contained in a charged 
partition.

\subsection{The bijection $\tau_l$, the ordering $\prec$ and abaci} \label{combinatorics_notation}

Throughout the proof of Theorem \ref{thm_A_J}, we need a large amount of notation which we introduce here. In particular, we have to pass from
$l$-multi-partitions (indexing the bases of the Fock space) to partitions (indexing the bases of the $q$-wedge space -- 
see Section \ref{section_fock_space}) and conversely. Following \cite{U2}, we achieve this using
a bijection $\tau_l$ which can be described in a combinatorial way (see Definition \ref{definition_tau_l}). This map is a variant of 
the bijection associating to a partition its $l$-quotient and its $l$-core. We construct here $\tau_l$ using abaci; for another (equivalent) 
description of $\tau_l$ and examples, see \cite[Remark 4.2 (ii) and Example 4.3]{U2}. The bijection $\tau_l$ is used in particular for defining 
the partial ordering $\prec$ on $\Pilm$ mentioned in Section \ref{section_our_results}; see Definition \ref{def_prec}.  

\subsubsection{Notation} \label{section_notation_tau_l}

The Euclidean algorithm shows that any integer $k \in \mZ$ can be written in a unique way as 
\begin{equation}
k=c(k)+n(d(k)-1)+nlm(k),
\end{equation}
with $c(k) \in \interv{1}{n}$, $d(k) \in \interv{1}{l}$ and $m(k) \in \mZ$. Consider the map
\begin{equation}
\phi : \mZ \rightarrow \mZ, \qquad k \mapsto c(k)+nm(k).
\end{equation}
$\phi$ enjoys the following obvious  properties, which we need later: for all $k, k' \in \mZ$, we have

\begin{eqnarray}
\label{property_P1}  \phi(k) \equiv c(k) \equiv k \pmod n, \\[5mm] 
\label{property_P2}  \bigl(k<k',\,  d(k)=d(k') \bigr) \Longrightarrow \phi(k)<\phi(k'),  \\[5mm]
\label{property_P3}  \bigl(k \leq k',\, \phi(k) \geq \phi(k') \bigr) \Longrightarrow m(k)=m(k').
\end{eqnarray}

\vspace{8mm}

\! \,For any $r$-tuple $\textbf{k}=(k_1,\ldots,k_r) \in \mZ^r$, let
\begin{equation} 
\textbf{c}(\textbf{k}):=(c(k_1),\ldots,c(k_r)) \in \mZ^r,
\end{equation}
and define in a similar way $\textbf{d}(\textbf{k})$. The group $\Sr$ acts on the left on $\mZ^{r}$ by
\begin{equation}
\sigma.(k_1,\ldots,k_r)=(k_{\sigma^{-1}(1)},\ldots,k_{\sigma^{-1}(r)}) \qquad ((k_1,\ldots,k_r) \in \mZ^r, \ \sigma \in \Sr),
\end{equation}
and a fundamental domain for this action is
$B := \set{(b_1, \ldots,b_r) \in \mZ^r \ | \ b_1 \geq \cdots \geq b_r}$. Let $\textbf{b}(\textbf{k})$ denote the element of $B$ that is conjugated 
to $\textbf{d}(\textbf{k})$ under the action of $\Sr$, $W_{\textbf{k}}$ be the stabilizer of $\textbf{b}(\textbf{k})$ (this is a parabolic subgroup 
of $\Sr$) and $\omega(\textbf{k})$ be the element of maximal length in $W_{\textbf{k}}$. Let $W^{\textbf{k}}$ be the set of minimal length 
representatives in the left cosets $\Sr/W_{\textbf{k}}$, and $v(\textbf{k})$ be the element in $W^{\textbf{k}}$ such that 
$\textbf{d}(\textbf{k})=v(\textbf{k}).\textbf{b}(\textbf{k})$. 

\begin{example}
Let $n=3$, $l=2$, $r=4$ and $\textbf{k}=(12, \! - \! 5,2,17)$. Then we have: \\

\noindent $\textbf{c}(\textbf{k})=(3,1,2,2),\ \textbf{d}(\textbf{k})=(2,1,1,2),\ \textbf{b}(\textbf{k})=(2,2,1,1),\ 
v(\textbf{k})=\sigma_3\sigma_2,\ \omega(\textbf{k})=\sigma_1\sigma_3.$ \fini
\end{example}

\begin{remark} \label{remark_permutation_v} 
Let $\textbf{k}=(k_1,\ldots,k_r) \in \mZ^r$. We can describe the action of $v(\textbf{k})^{-1}$ on $\textbf{k}$ as follows. Consider 
$\textbf{k}$ as a word formed by the letters $k_i$ and for $1 \leq d \leq l$, denote by $w_d$ the subword of $\textbf{k}$ formed by the 
letters $k_i$ such that $d(k_i)=d$. Then we have $v(\textbf{k})^{-1}.\textbf{k}=w_l \cdots w_1$.~\fini
\end{remark}

\subsubsection{The bijection $\tau_l$, the ordering $\prec$ and abaci}

\begin{definition} 
A \emph{$1$-runner abacus} is a subset $A$ of $\mZ$ such that $-k \in A$ and $k \notin A$ for all large enough $k \in \mN$. In a less formal way, 
each $k \in A$ corresponds to the position of a bead on the horizontal abacus $A$ which is full of beads on the left and empty on the right. 
Let $\mathcal{A}$ be the set of $1$-runner abaci. If $N \geq 1$, an \emph{$N$-runner abacus} is an $N$-tuple of $1$-runner abaci. If
$\textbf{A}=(A_1,\ldots,A_N) \in \mathcal{A}^N$ is an $N$-runner abacus, we identify $\textbf{A}$ with the subset
\begin{equation}
\set{  (k,d) \ | \ 1 \leq d \leq N, \ k \in A_d} \subset \mZ \times \interv{1}{N}.
\end{equation}
\fini
\end{definition}

To $\boldsymbol{\lambda}=(\lambda^{(1)},\ldots,\lambda^{(N)}) \in \Pi^N$ and $\textbf{s}=(s_1,\ldots,s_N) \in \mZ^N$ we associate the $N$-runner 
abacus
\begin{equation}
A(\boldsymbol{\lambda}, \textbf{s}):=\set{(\lambda_i^{(d)}+s_d+1-i, \, d) \ | \ i \geq 1, \ 1 \leq d \leq N}.
\end{equation}
One checks easily that the map
\begin{equation}
(\boldsymbol{\lambda}, \textbf{s}) \in \Pi^N \times \mZ^N \mapsto A(\boldsymbol{\lambda}, \textbf{s}) \in \mathcal{A}^N
\end{equation}
is bijective. \\

Recall the definition of the maps $k \mapsto d(k)$ and $k \mapsto \phi(k)$ from Section \ref{section_notation_tau_l}.
Note that $k \in \mZ \mapsto \big( \phi(k),d(k) \big) \in \mZ \times \interv{1}{l}$ is a bijection.

\begin{definition} \label{definition_tau_l}
The bijection $\tau_l : \Pi \times \mZ \cong \mathcal{A} \rightarrow \Pi^l \times \mZ^l \cong \mathcal{A}^l$ is defined in 
terms of abaci by the formula
\begin{equation}
\tau_l(A):= \big\{\big(\phi(k),d(k)\big) \ \big| \ k \in A \big\} \in \mathcal{A}^l \qquad (A \in \mathcal{A}). 
\end{equation} \fini
\end{definition}

\begin{remark} Let $\lambdal \in \Pi^l$, $\textbf{s}_l=(s_1,\ldots,s_l) \in \mZ^l$, $\lambda \in \Pi$ and $s \in \mZ$ satisfying the relation
$(\lambdal,\textbf{s}_l)=\tau_l(\lambda,s).$ Then we have $s=s_1+\cdots+s_l$. \fini
\end{remark}

\begin{notation} \label{notation_lambda_corresponding_lambdal}
Let $\esselle=(s_1,\ldots,s_l) \in \mZ^l$ and $s:=s_1+\cdots+s_l$. Write
\begin{equation} 
\lambdal \stackrel{\esselle}{\longleftrightarrow} \lambda
\end{equation}
if $\lambda \in \Pi$ and $\lambdal \in \Pilm$ are related by $(\lambdal,\esselle) = \tau_l(\lambda,s)$. We drop the $\esselle$ in the
notation if it is clearly given by the context. \fini 
\end{notation}

\begin{example} 
Let $n=2$, $l=3$, $m=5$ and $\esselle=(0,0,-1)$. Then Figure \ref{fig:1} shows that 
$$\bigl((1,1),(1,1),(1) \bigr) \stackrel{\esselle}{\longleftrightarrow} (4,3,3,2,1).$$ 

%
\vspace*{5mm} 
\begin{figure}[htbp]
\begin{center}
\includegraphics[height=6.5cm]{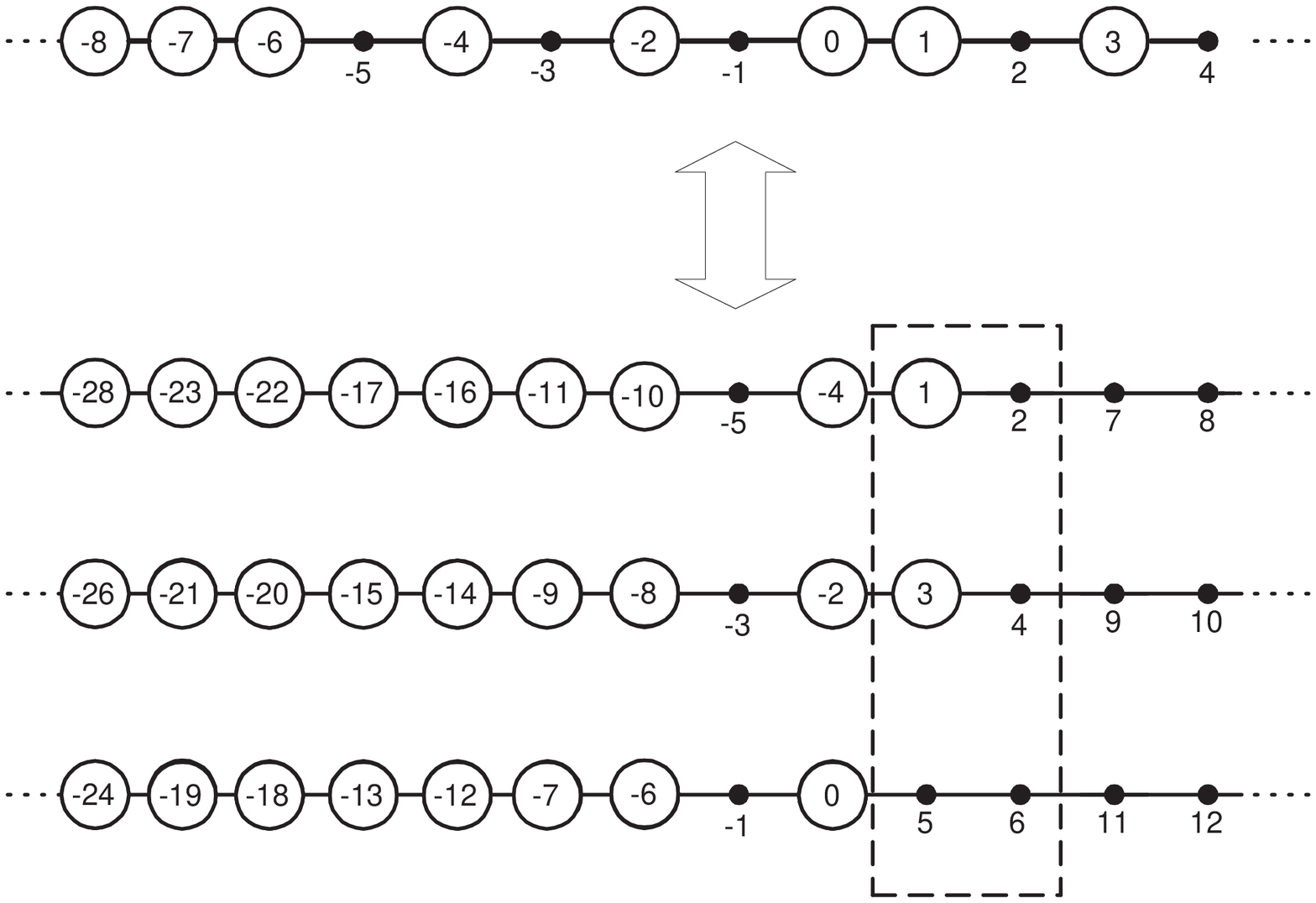} 
\end{center}
\caption{Computation of the bijection $\tau_l$ using abaci.}
\label{fig:1} 
\end{figure}
\fini
\end{example}

We now define a partial ordering $\prec$ on $\Pilm$ as follows. 

\begin{definition} \label{def_prec}
Let $\esselle=(s_1,\ldots,s_l) \in \mZ^l$. Let $\lambdal$, $\mul \in \Pilm$ and $\lambda,  \mu \in \Pi$ 
be such that $\lambdal \stackrel{\esselle}{\longleftrightarrow} \lambda$ and 
$\mul \stackrel{\esselle}{\longleftrightarrow} \mu$. We say that $\lambdal$ \emph{precedes} $\mul$ and write
\begin{equation}
\lambdal \preceq \mul
\end{equation}
if $\mu$ dominates $\lambda$. In particular, by (\ref{condition1_dominance_ordering}), $\lambda$ and $\mu$ must be 
partitions of the same integer. Note that the ordering $\preceq$ depends on the multi-charge $\esselle$ that we consider. 
Write $\lambdal \prec \mul$ if $\lambdal \preceq \mul$ 
and $\lambdal \neq \mul$. \fini
\end{definition} 

\subsection{$\beta$-numbers and ribbons} 

Throughout this section we fix an integer $s \in \mZ$.

\begin{definition} \label{def_beta_numbers}
Let $\lambda=(\lambda_1,\lambda_2,\ldots) \in \Pi$ be a partition
with at most $r$ parts. The $r$-tuple 
\begin{equation}
\boldsymbol{\beta}_r(\lambda):=(\lambda_1+s,\lambda_2+s -1,\ldots,\lambda_r+s-r+1) \in \mZ^r
\end{equation}
is called the \emph{$r-$list of $\beta$-numbers associated to} $(\lambda,s)$ or (with a slight abuse of notation) the \emph{list} or 
\emph{sequence of $\beta$-numbers associated to $\lambda$}. The set of integers that form $\boldsymbol{\beta}_r(\lambda)$ is denoted by 
$B_r(\lambda)$. \fini
\end{definition}

With the notation of the definition above, note that $\boldsymbol{\beta}_r(\lambda)$ is a decreasing sequence of integers all greater than (or 
equal to) $s+1-r$. This sequence depends on the integer $s$ we have fixed, but we do not mention it in our notation. Note that a partition 
$\lambda$ is completely determined by its sequence of $\beta$-numbers. If $r=|\lambda|$, write more simply 
\begin{equation}
\boldsymbol{\beta}(\lambda):=\boldsymbol{\beta}_r(\lambda) \qquad \mbox{and} \qquad B(\lambda):=B_r(\lambda).
\end{equation}
If $f$ is a function defined on $\mZ^r$, it is convenient to consider $f$ as a function (still denoted by $f$) defined
on the set of partitions of $r$ by the formula 
\begin{equation}
f(\lambda) :=f (\boldsymbol{\beta}(\lambda)) \qquad (\lambda \in \Pi, \,  |\lambda|=r).
\end{equation} 
\! \,For example, we define this way for any partition $\lambda$ the vectors $\textbf{c}(\lambda)$, $\textbf{d}(\lambda)$ and so on. 
See Section \ref{section_notation_tau_l} for the corresponding notation. \\ 

In order to prove Theorem \ref{thm_A_J}, we have to relate the adding/removal of a ribbon in a charged partition and the corresponding
$\beta$-numbers. Let us recall a classical result on $\beta$-numbers (see \eg \cite[Lemma 5.26]{Mat1}). 

\begin{lemma} \label{ribbon_beta_numbers_1} Let $\nu$ and $\kappa$ be two partitions with at most $r$ parts, and let
$\boldsymbol{\beta}_r(\nu)=(\alpha_1,\ldots,\alpha_r)$ and $\boldsymbol{\beta}_r(\kappa)=(\beta_1,\ldots,\beta_r)$ denote the sequences of 
$\beta$-numbers associated to $\nu$ and $\kappa$ respectively. Then the following statements are equivalent.
\begin{itemize}
\item[\rm (i)] $\nu \subset \kappa$, and $\rho:=\kappa / \nu$ is a ribbon of length $h$.
\item[\rm (ii)] There exist positive integers $b$ and $h$ such that 
$B_r(\nu)=\set{\beta_1,\ldots,\beta_{b-1},\beta_b-h,\beta_{b+1},\ldots,\beta_r}.$
\end{itemize} 
In this case, $b$ is the row number of the tail of $\rho$ and $h$ is the length of $\rho$. Let $\sigma \in \Sr$ denote the permutation obtained by
arranging decreasingly the integers $(\beta_1,\ldots,\beta_{b-1},\beta_b-h,\beta_{b+1},\ldots,\beta_r).$ Then we have 
$\ell(\sigma)=\haut(\rho)$. Moreover, the content of the head of $\rho$ is 
\begin{equation}
\cont(\hd(\rho))=\alpha_c = \beta_b-h,
\end{equation}
where $c$ is the row number of the head of $\rho$.
\end{lemma}

\begin{demo} 
The proof of (i) $\Rightarrow$ (ii) is easy. Conversely, assume that (ii) holds. Then we must have $\beta_b-h \geq s+1-r$, and there must exist
$b \leq c \leq r$ such that $\beta_c > \beta_b-h > \beta_{c+1}$ (if $c=r$, put $\beta_{c+1}:=s-r$). Note then that $\nu$ is obtained from 
$\kappa$ by removing a ribbon $\rho$, where $\rho \subset \kappa$ is the ribbon whose head is located at row $c$ of $\kappa$ and whose tail is 
located at row $b$ of $\kappa$. $\rho$ is actually a ribbon of length $h$. Moreover, with the notation of the statement of this lemma, we have 
$\sigma.(\beta_1,\ldots,\beta_{b-1},\beta_b-h,\beta_{b+1},\ldots,\beta_r)
=(\beta_1,\ldots,\beta_{b-1},\beta_{b+1},\ldots,\beta_c,\beta_b-h,\beta_{c+1},\ldots,\beta_r),$
hence $\sigma$ is a cycle of length $c-b=\haut(\rho)$. Finally, the head of $\rho$ has coordinates $(c, \nu_c+1)$,
so its content is equal to $\cont(\hd(\rho)) = s + (\nu_c+1) - c = \alpha_c =\beta_b-h$.
\end{demo}

\begin{example}
Let $s=4$, $r=5$, $\kappa=(6,5,3,2,2)$ and $\nu=(6,2,2,2,2)$. Then the skew diagram $\rho := \kappa / \nu$ is a ribbon and we have 
$b=2$, $c=3$ and $h=4$. Moreover, we have $\boldsymbol{\beta}(\kappa)=(\beta_1,\ldots,\beta_5)=(10,8,5,3,2)$ and 
$\boldsymbol{\beta}(\nu)=(\beta_1,\beta_3,\beta_b-h,\beta_4,\beta_5)=(10,5,4,3,2)$. We have $\sigma=(2,3)$, hence $\ell(\sigma)=1=\haut(\rho)$. 
The head of $\rho$ has coordinates $(3,3)$, so its content is $\cont(\hd(\rho)) = 4 = \beta_b-h$. \fini
\end{example}

\begin{lemma} \label{ribbon_beta_numbers_2} Let $\nu$, $\kappa \in \Pi$ be such that $|\nu|=|\kappa|=r$ and $\nu \neq \kappa$. Let 
$\boldsymbol{\beta}(\nu)=(\alpha_1,\ldots,\alpha_r)$ and $\boldsymbol{\beta}(\kappa)=(\beta_1,\ldots,\beta_r)$ denote the sequences of 
$\beta$-numbers associated to $\nu$ and $\kappa$ respectively. Set $\rho:=\nu / (\nu \cap \kappa)$ and $\rho':=\kappa / (\nu \cap \kappa)$. \\ 
\begin{itemize}
\item[1)] Then, $\rho$ and $\rho'$ are two ribbons if and only if  $\sharp \bigl(B(\nu) \cap B(\kappa) \bigr) = r-2$. In this case, denote by
\begin{itemize}
\item[.] $h$ the common length of $\rho$ and $\rho'$,
\item[.] $y$ the row number of the tail of $\rho'$,
\item[.] $y'$ the row number of the head of $\rho'$,
\item[.] $x'$ the row number of the tail of $\rho$, and
\item[.] $x$ the row number of the head of $\rho$.
\end{itemize} \medskip

Then we have 
\begin{eqnarray}
\set{\alpha_i \mid i \neq x',y'} &=& \set{\beta_j \mid j \neq x,y}, \\
\cont(\hd(\rho))=\beta_x=\alpha_{x'}-h \quad & \mbox{and} & \quad \cont(\hd(\rho'))=\alpha_{y'} = \beta_y-h.
\end{eqnarray}
Let $\pi \in \Sr$ be the permutation obtained by arranging decreasingly the integers forming $B(\nu)$. Then we have
$\ell(\pi)=\haut(\rho)+\haut(\rho')$. \bigskip

\item[2)] Assume that the conditions of 1) hold. Then we have the following equivalences, and moreover one of the two following cases occurs:

$$\begin{array}{rccl}
\mbox{\rm (i)} &  y \leq y' < x' \leq x & \Longleftrightarrow & \nu \lhd \kappa  , \\
\mbox{\rm (ii)} & x' \leq x < y \leq y' & \Longleftrightarrow & \kappa \lhd \nu.
\end{array}$$

\end{itemize}
\end{lemma}

\begin{demo} 
We prove 1) by applying the previous lemma to the pairs of partitions $(\nu \cap \kappa,\nu)$ and $(\nu \cap \kappa,\kappa)$. Let us prove 2). 
The inequalities $y \leq y'$ and $x' \leq x$ are obvious. Since $\rho \cap \rho' = \emptyset$, one of the two following cases occurs: either 
$y'<x'$ and then $\nu \lhd \kappa$, or $x<y$ and then $\kappa \lhd \nu$. This proves both implications $\Rightarrow$, and since one of the two 
cases occurs, we get the desired equivalences.
\end{demo}

\subsection{Definition of the matrix $J^{<}$} \label{def_J} 

Let $R$ be a local ring, with unique maximal ideal $\wp$. We define in this section a matrix
$\mathcal{J}=\bigl(J_{\lambdal,\mul} \bigr)_{\lambdal, \mul \in \Pilm}$, with coefficients in $R$, depending on parameters
$m$, $l \in \mN^*$ and $v$, $u_1,\ldots,u_l \in R$. This matrix is closely related to the matrix formed by the entries 
$\nu_{\wp}(g_{\lambdal,\mul})$ of \cite{JM} (see (\ref{eq_matrix_J_Jantzen_sum_formula})).
Let $\lambdal=(\lambda^{(1)},\ldots,\lambda^{(l)})$, $\mul = (\mu^{(1)},\ldots,\mu^{(l)}) \in \Pilm$, and consider the following cases.

\begin{figure}[htbp]
\begin{center}
\includegraphics[height=4cm]{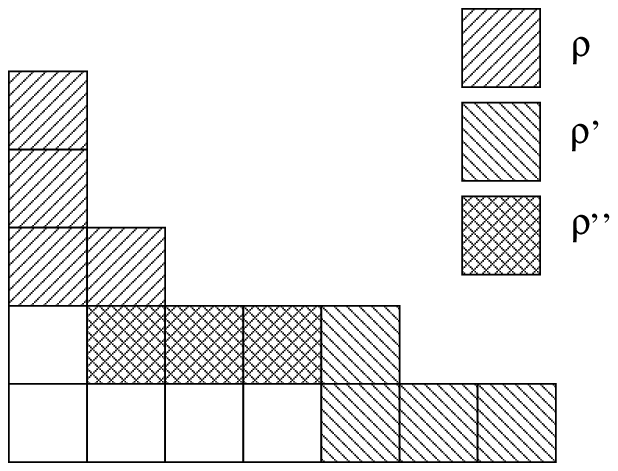}
\end{center}
\caption{The ribbons $\rho$, $\rho'$ and $\rho''$ (the nodes of $(\lambda^{(d)} \cap \mu^{(d)}) - \rho''$ are depicted in white).} 
\label{fig:2} 
\end{figure}

\begin{itemize} 
\item \underline{\emph{Case $(J_1)$.}} Assume that $\lambdal \neq \mul$ and that there exist two integers $d$, $d' \in \interv{1}{l}$, $d \neq d'$
satisfying the following conditions: $\mu^{(d)} \subset \lambda^{(d)}$, $\lambda^{(d')} \subset \mu^{(d')}$, $\lambda^{(b)}=\mu^{(b)}$ for all 
integer $b \in \interv{1}{l} \setminus \set{d,d'}$, and $\rho:=\lambda^{(d)} / \mu^{(d)}$ and $\rho':=\mu^{(d')} / \lambda^{(d')}$ are two ribbons 
of the same length $\widehat{h}$. Let $\hd(\rho)=(i,j,d)$ denote the head of $\rho$ and $\hd(\rho')=(i',j',d')$ denote the head of $\rho'$. Set
\begin{equation}
\varepsilon:=(-1)^{\haut(\rho)+\haut(\rho')} \qquad \mbox{and}
\end{equation}
\begin{equation} 
J_{\lambdal,\mul}:={\bigl(u_d v^{j-i}-u_{d'}v^{j'-i'}\bigr)}^{\varepsilon}.
\end{equation}
  
\item \underline{\emph{Case $(J_2)$.}} Assume that $\lambdal \neq \mul$ and that there exists $d \in \interv{1}{l}$ such that
$\lambda^{(b)}=\mu^{(b)}$ for all $b \neq d$, and $\rho:=\lambda^{(d)} / (\lambda^{(d)} \cap \mu^{(d)})$ and 
$\rho':=\mu^{(d)} / (\lambda^{(d)} \cap \mu^{(d)})$ are two ribbons of the same length $\widehat{h}$.
By definition of $\rho$ and $\rho'$, we have $\rho \cap \rho' = \emptyset$, whence we get (depending on the relative positions of $\rho$ and 
$\rho'$) that either $\lambda^{(d)} \lhd \mu^{(d)}$ or $\mu^{(d)} \lhd \lambda^{(d)}$. Assume that $\lambda^{(d)} \lhd \mu^{(d)}$. Let 
$\rho'' \subset (\lambda^{(d)} \cap \mu^{(d)})$ be the ribbon obtained by connecting the tail of $\rho$ to the head of $\rho'$, excluding the two 
latter nodes (see Figure \ref{fig:2}). Denote by $\hd(\rho)=(i,j,d)$ (\resp $\hd(\rho')=(i',j',d')$, \resp $\hd(\rho'')=(i'',j'',d'')$) the head of 
$\rho$ (\resp $\rho'$, \resp $\rho''$), and finally set  

\begin{equation}
\varepsilon_1:=(-1)^{\haut(\rho)+\haut(\rho')}, \qquad \varepsilon_2:=(-1)^{\haut(\rho \cup \rho'')+\haut(\rho'' \cup \rho')} \qquad \mbox{and}
\end{equation}
\begin{equation}
J_{\lambdal,\mul}:= {\bigl(u_d (v^{j-i} - v^{j'-i'})\bigr)}^{\varepsilon_1}.{\bigl(u_d (v^{j-i} - v^{j''-i''})\bigr)}^{\varepsilon_2}.
\end{equation}

If $\mu^{(d)} \lhd \lambda^{(d)}$, set $J_{\lambdal,\mul}:=J_{\mul,\lambdal}$. 

\medskip

\item \underline{\emph{Case $(J_3)$.}} In all other cases, set 
\begin{equation}
J_{\lambdal,\mul}:=1.
\end{equation}
\end{itemize} 

%
\vspace*{1cm} 

\vspace{5mm}
We now define a matrix $J=J_{\wp}=\bigl(j_{\lambdal,\mul} \bigr)_{\lambdal, \mul \in \Pilm}$, with integer 
coefficients, by the formula    
\begin{equation}
j_{\lambdal,\mul} := \nu_{\wp}(J_{\lambdal,\mul}) \qquad (\lambdal, \mul \in \Pilm).
\end{equation}

\noindent Now, let $\leqslant$ be an arbitrary partial ordering on $\Pilm$ and write $\lambdal < \mul$ if $\lambdal \leqslant \mul$ and 
$\lambdal \neq \mul$ ($\lambdal,\, \mul \in \Pilm$). Recall the definition of the matrix 
$J^{<}=\bigl(j^{<}_{\lambdal,\mul} \bigr)_{\lambdal,\, \mul \in \Pilm}$ from (\ref{eq_def_J_inf})\,; namely, put
\begin{equation} 
j^{<}_{\lambdal,\mul}:= \left\{ 
\begin{array}{cc}
j_{\lambdal,\mul} & \mbox{if } \lambdal < \mul \\[1mm]
0 & \mbox{otherwise}
\end{array}
\right. \qquad (\lambdal,\, \mul \in \Pilm).
\end{equation}

If we take $\leqslant \ = \ \unlhd$, then we get a matrix $J^{\lhd}$ whose entries correspond, up to conjugation of multi-partitions,
to the integers $\nu_\wp(g_{\lambdal,\mul})$ of \cite{JM} (see (\ref{eq_matrix_J_Jantzen_sum_formula})). Given a multi-charge $\esselle$, we shall 
also consider the matrix $J^{\prec}$, where the ordering $\prec$ (depending on $\esselle$) was introduced in Definition \ref{def_prec}. This is 
the matrix $J^{\prec}$ of Theorems \ref{thm_Delta_J} and \ref{thm_A_J}. If $\esselle$ is $m$-dominant (in the sense of Definition
\ref{def_dominant}), then the matrices $J^{\prec}$ and $J^{\lhd}$ coincide (see Proposition \ref{expr_J_dominant}).

\section{$q$-deformed higher-level Fock spaces} 

In this section we follow \cite{U2}, to which we refer the reader for more details. The vector spaces we consider here are over 
$\mC(q)$, where $q$ is an indeterminate over $\mC$.  
 
\subsection{$q$-wedge products and higher-level Fock spaces} \label{section_fock_space}

Let $s \in \mZ$. Let $\Lambda^s$ denote the (semi-infinite) $q$-wedge space of charge $s$ (this space is denoted by $\Lambda^{s+\frac{\infty}{2}}$ 
in \cite{U2}). $\Lambda^s$ is an integrable representation of level $l$ of the quantum algebra $\Uq$. As a vector space, it has a natural basis 
formed by the so-called \emph{ordered $q$-wedge products}. These vectors can be written as 
\begin{equation}
u_{\textbf{k}} = u_{k_1} \wedge u_{k_2} \wedge \cdots,
\end{equation}
where $\textbf{k}=(k_i)_{i \geq 1}$ is a decreasing sequence of integers such that $k_i = s + 1 - i$ for $i \gg 0$. The basis formed by the 
ordered wedge products is called \emph{standard}. More generally, we use the 
\emph{non-ordered wedge products}; a non-ordered wedge product $u_{\textbf{k}} = u_{k_1} \wedge u_{k_2} \wedge \cdots \ \in \Lambda^s$ is indexed 
by a sequence of integers $(k_i)$ such that $k_i = s + 1 - i$ for $i \gg 0$, but we no longer require that $(k_i)$ is decreasing. Any non-ordered 
wedge product can be written as a linear combination of ordered wedge products by using the so-called \emph{ordering rules}, which are given in 
\cite[Proposition 3.16]{U2} and in a slightly different form in Proposition \ref{base_v_k}. \\

The vectors of the standard basis of $\Lambda^s$ can also be indexed by partitions as follows. Let 
$u_{\textbf{k}} = u_{k_1} \wedge u_{k_2} \wedge \cdots \ \in \Lambda^s$ be an ordered wedge product. For $i \geq 1$ set 
$\lambda_i:=k_i-s+i-1$; then $\lambda:=(\lambda_1,\lambda_2,\ldots)$ is a partition. We then write $u_{\textbf{k}} = |\lambda,s \rangle$.
Note that if $\lambda$ has at most $r$ parts, then we have $(k_1,\ldots,k_r)=\boldsymbol{\beta}_r(\lambda)$, which explains the definition of
the $\beta$-numbers we gave in Definition \ref{def_beta_numbers}. \\

Let $\Fq$ be the higher-level Fock space with multi-charge $\esselle=(s_1,\ldots,s_l) \in \mZ^l$ \cite{U2}. As a vector space, $\Fq$ has a 
natural basis $\set{ |\lambdal,\esselle \rangle \mid \lambdal \in \Pi^l}$ indexed by $l$-multi-partitions. If $s=s_1+\cdots+ s_l$, then $\Fq$
can be identified with a subspace of $\Lambda^s$ by the embedding $\Fq \hookrightarrow \Lambda^s$, 
$|\lambdal,\esselle \rangle \mapsto |\lambda,s \rangle$, where $\lambda$ is the partition such that $\lambdal \leftrightarrow \lambda$ (see 
Notation \ref{notation_lambda_corresponding_lambdal} for the meaning of $\leftrightarrow$). We make from now on this identification; in fact, 
$\Lambda^s$ is isomorphic to the direct sum of all the $\textbf{F}_q[\textbf{t}_l]$'s, where $\textbf{t}_l$ is any $l$-tuple of 
integers summing to $s$. Thus, the vectors of the standard basis of $\Lambda^s$ can also be indexed by charged $l$-multi-partitions.

\subsection{The involution $\barre$} \label{section_fock_involution}

In order to define the canonical basis of $\Lambda^s$, we equip this space with an involution $\barre$.

\begin{definition} The involution $\barre$ of  $\Lambda^s$ is the $\mC$-vector
space automorphism that maps $q$ to $q^{-1}$ and that acts on the standard basis of $\Lambda^s$ as follows 
\cite[Proposition 3.23 and Remark 3.24]{U2}. Let $\lambda \in \Pi$ be a partition of $r$, and $\textbf{k}=(k_i) \in \mZ^{\mN^{*}}$ be such that   
$u_{\textbf{k}}=|\lambda,s \rangle$. Then 

\begin{equation} \label{eq_action_barre}
\overline{|\lambda,s \rangle}:=(-1)^{\kappa(\textbf{d}(\lambda))}q^{\kappa(\textbf{d}(\lambda))-\kappa(\textbf{c}(\lambda))} 
(u_{k_r} \wedge \cdots \wedge u_{k_1}) \wedge u_{k_{r+1}} \wedge u_{k_{r+2}} \wedge \cdots,
\end{equation}

\bigskip

\noindent where for any $\textbf{a}=(a_1,\ldots,a_r) \in \mZ^r$, $\kappa(\textbf{a})$ is the integer defined by
\begin{equation}
\kappa(\textbf{a}) := \sharp \set{(i,j) \in \mN^2 \mid 1 \leq i < j \leq r, \ a_i=a_j},
\end{equation} 
and $\textbf{c}(\lambda)$ and $\textbf{d}(\lambda)$ are defined in Section \ref{section_notation_tau_l}. \fini
\end{definition}

We can straighten the non-ordered wedge product in the right-hand side of (\ref{eq_action_barre}) in order to 
express it as a linear combination of ordered wedge products. \\

One checks that $\barre$ preserves the subspace 
\begin{equation}
\Fq_m :=  \bigoplus_{\lambdal \in \Pilm}{\mC(q)\, |\lambdal, \esselle \rangle} \subset \Fq. 
\end{equation} 

\begin{definition}
Define a matrix $A(q)=\bigl(a_{\lambdal,\mul}(q) \bigr)_{\lambdal, \mul \in \Pilm}$ with entries in $\mC(q)$ by

\begin{equation}
\overline{|\mul, \esselle \rangle} = 
\sum_{\lambdal \in \Pilm}{a_{\lambdal,\mul}(q)\, |\lambdal, \esselle \rangle} \qquad (\mul \in \Pilm).   
\end{equation} \fini
\end{definition}

\noindent The matrix $A(q)$ depends on $n$, $l$, $\esselle$ and $m$. The ordering rules show that $A(q)$ is unitriangular
with respect to $\preceq$, that is 
\begin{equation} \label{eq_barre_unitr}
a_{\lambdal,\mul}(q) \neq 0 \Rightarrow \lambdal \preceq \mul \qquad \mbox{and} \qquad a_{\lambdal,\lambdal}(q) =1 \qquad 
(\lambdal, \mul \in \Pilm). 
\end{equation}

\noindent The same rules also imply that $A(1)$ is the identity matrix.      

\subsection{Uglov's canonical basis} \label{section_canonical_basis}

Since the matrix $A(q)$ of the involution of $\Fq_{m}$ is unitriangular, a classical argument can be used to prove the following result.

\begin{thm}[\cite{U2}] \label{thm_bc}
There exists a unique basis $\set{\mathcal{G}(\lambdal,\esselle) \mid \lambdal \in \Pilm}$ of $\Fq_{m}$ satisfying both 
following conditions:

\begin{itemize} 
\item[\rm (i)] $\overline{\mathcal{G}(\lambdal,\esselle)}=\mathcal{G}(\lambdal,\esselle)$,
\item[\rm (ii)] $\ds \mathcal{G}(\lambdal,\esselle) - |\lambdal,\esselle \rangle \ \in \ \bigoplus_{\mul \in \Pilm}{q \, \mC[q] \,
|\mul,\esselle \rangle}.$ \cqfd
\end{itemize}
\end{thm}

\begin{definition}
The basis $\set{\mathcal{G}(\lambdal,\esselle) \mid \lambdal \in \Pilm}$ is called the \emph{canonical basis} of $\Fq_{m}$. Define a matrix
$\Delta(q)=\bigl(\Delta_{\lambdal,\mul}(q) \bigr)_ {\lambdal, \mul \in \Pilm}$ with entries in $\mC[q]$ by 
\begin{equation}
\mathcal{G}(\mul,\esselle) = 
\sum_{\lambdal \in \Pilm}{\Delta_{\lambdal,\mul}(q)\, |\lambdal, \esselle \rangle} \qquad (\mul \in \Pilm).
\end{equation} \fini
\end{definition}

\noindent The matrix $\Delta(q)$ depends on $n$, $l$, $\esselle$ and $m$. By Condition (ii) of Theorem \ref{thm_bc}, the matrix 
$\Delta(q)$ is also unitriangular with respect to $\preceq$. By \cite[Theorem 3.26]{U2}, the entries of $\Delta(q)$ can be expressed as 
Kazhdan-Lusztig polynomials related to parabolic modules of an affine Hecke algebra of type $\tilde{A}$, so by \cite{KT}, these entries are in 
$\mN[q]$.

\subsection{Another basis of $\Lambda^s$. Ordering rules.} \label{base_v_k}

The ordering rules $(R_1)$-$(R_4)$ from \cite[Proposition 3.16]{U2} do not give at $q=1$ anticommuting relations like
$u_{k_1} \wedge u_{k_2} = - u_{k_2} \wedge u_{k_1}$, because of the signs involved in Rules $(R_3)$ and $(R_4)$. To fix this, we introduce another
basis of $\Lambda^s$ that differs from the standard basis only by signs. The basis we consider here is actually the basis of ordered
wedge products introduced in \cite{U1}. $\Lambda^s$ is graded by
\begin{equation}
\deg(|\lambda,s \rangle):=|\lambda| \qquad (\lambda \in \Pi). 
\end{equation}  

\bigskip 

\noindent Let $u_{\textbf{k}} =  u_{k_1} \wedge u_{k_2} \wedge \cdots \in \Lambda^s$ be a (not necessarily ordered) wedge product of degree $r$. 
Set

\begin{equation}
\begin{array}{crl}
& v_{\textbf{k}} = v_{k_1} \wedge v_{k_2} \wedge \cdots & := (-1)^{\ell(v(k_1,\ldots,k_r))} u_{\textbf{k}} \\[2mm]
\mbox{and similarly} & v_{k_1} \wedge \cdots \wedge v_{k_r} &:= (-1)^{\ell(v(k_1,\ldots,k_r))} u_{k_1} \wedge \cdots \wedge u_{k_r},
\end{array}
\end{equation}
where $v(k_1,\ldots,k_r) \in \Sr$ is defined in Section \ref{section_notation_tau_l}. (If $\textbf{k}=(k_1,\ldots,k_r) \in \mZ^r$, we hope that the 
reader will make easily the difference between the permutation $v(\textbf{k}) \in \Sr$ and the wedge product 
$v_{\textbf{k}}=v_{k_1} \wedge \cdots \wedge v_{k_r}$.) We say that the wedge product $v_{\textbf{k}}$ is \emph{ordered} if so is $u_{\textbf{k}}$. 
It is straightforward to see, using the ordering rules for the $u_{\textbf{k}}$'s given by \cite[Proposition 3.16]{U2}, that the ordering rules for 
the $v_{\textbf{k}}$'s are given by the following proposition.

\begin{prop} \ \\ 
\rm (i) \it Let $k_1 \leq k_2$, and $\gamma \in \interv{0}{nl-1} \ ( \mbox{\resp} \delta \in \interv{0}{nl-1})$ denote the residue of 
$c(k_2)-c(k_1)$ $\bigl( \bigr.$\resp of $n(d(k_2)-d(k_1)) \bigl. \bigr)$ modulo $nl$. Then we have
\begin{equation}
v_{k_1} \wedge v_{k_2} = - v_{k_2} \wedge v_{k_1} \hspace{7cm} \mbox{if } \gamma=\delta=0, 
\tag{$R_1$} \end{equation}  
\begin{equation}
 \begin{array}{rll} v_{k_1} \wedge v_{k_2} &=& - q^{-1} v_{k_2} \wedge v_{k_1} \\[5mm]
&&- (q^{-2}-1) \ds \sum_{i \geq 1}{q^{-2i+1} v_{k_2-nli} \wedge v_{k_1+nli}} \\[5mm]
&& + (q^{-2}-1) \ds \sum_{i \geq 0}{q^{-2i} v_{k_2-\gamma-nli} \wedge v_{k_1+\gamma+nli}} \\
\end{array} \\
\qquad \qquad \mbox{if } \gamma>0,  \delta=0,
\tag{$R_2$} \end{equation}
\begin{equation}
\begin{array}{rll} v_{k_1} \wedge v_{k_2} &=& -q v_{k_2} \wedge v_{k_1} \\[5mm]
&&- (q^2-1) \ds \sum_{i \geq 1}{q^{2i-1} v_{k_2-nli} \wedge v_{k_1+nli}} \\[5mm]
&& + (q^2-1) \ds \sum_{i \geq 0}{q^{2i} v_{k_2-\delta-nli} \wedge v_{k_1+\delta+nli}} \\
\end{array} \\
\qquad \qquad \qquad \mbox{if } \gamma=0,  \delta>0,
\tag{$R_3$} \end{equation}
\begin{equation}
\begin{array}{rll} v_{k_1} \wedge v_{k_2} &=& -v_{k_2} \wedge v_{k_1} \\[5mm]
&&- (q-q^{-1}) \ds \sum_{i \geq 1}{\frac{q^{2i}-q^{-2i}}{q+q^{-1}}  v_{k_2-nli} \wedge v_{k_1+nli}}  \\[5mm]
&&-(q-q^{-1}) \ds \sum_{i \geq 0}{\frac{q^{2i+1}+q^{-2i-1}}{q+q^{-1}}  v_{k_2-\gamma-nli} \wedge v_{k_1+\gamma+nli}} \\[5mm]
&&+(q-q^{-1}) \ds \sum_{i \geq 0}{\frac{q^{2i+1}+q^{-2i-1}}{q+q^{-1}}  v_{k_2-\delta-nli} \wedge v_{k_1+\delta+nli}} \\[5mm]
&&+(q-q^{-1}) \ds \sum_{i \geq 0}{\frac{q^{2i+2}-q^{-2i-2}}{q+q^{-1}}  v_{k_2-\gamma-\delta-nli} \wedge v_{k_1+\gamma+\delta+nli}} \\[5mm]
\end{array}
\quad \mbox{if } \gamma>0,  \delta>0, 
\tag{$R_4$} \end{equation} \\

\noindent where the sums range over the indices $i$ such that the corresponding wedge products are ordered. \\

\noindent \rm (ii) \it The rules from \rm (i) \it are valid for any pair of adjacent factors of the $q$-wedge product 
$v_{\textbf{k}} = v_{k_1} \wedge v_{k_2} \cdots$. \cqfd
\end{prop}

\vspace{5mm}

Let us end this section by a useful piece of notation.
\bigskip 

\begin{notation} \label{notation_permuted_wedges}
Let $\nu \in \Pi$ be a partition of $r$ and $\sigma \in \Sr$. Set
\begin{equation}
u_{\sigma.\nu}:=u_{\sigma.\boldsymbol{\beta}(\nu)} \qquad \mbox{and similarly} \qquad  v_{\sigma.\nu}:=v_{\sigma.\boldsymbol{\beta}(\nu)}
\end{equation}
(these are wedge products of $r$ factors each). We say that $u_{\sigma.\nu}$ (\resp $v_{\sigma.\nu}$) is \emph{obtained from $u_{\nu}$
(\resp $v_{\nu}$) by permutation}.  
 \fini
\end{notation}

\vspace{5mm}
\begin{center}
\scshape \Large PART C: Proof of Theorem \ref{thm_A_J}  
\end{center}
\vspace{5mm}

We now start the proof of Theorem \ref{thm_A_J}. In Section \ref{section_expression_J}, we give a simpler expression for the entries of the matrix
$J^{\prec}$ (see Proposition \ref{proposition_J1_J2}). 
In Section \ref{section_good_sequence}, we compute the derivative at $q=1$ of the $\barre$ involution of $\Fq_m$ in terms of good sequences that
we introduce in Definition \ref{def_good_sequence}; the result is given in Proposition \ref{express_aprime_lambda_mu}. We compare both expressions 
in Section \ref{proof_thm_A_J} in order to complete the proof. Apart from this, we compare in Section \ref{section_case_sl_dominant} the matrices 
$J^{\prec}$ and $J^{\lhd}$ when the multi-charge $\esselle$ is $m$-dominant. \\

\noindent \textbf{Notation for Part C}\! \,From now on, we consider the modular system $(R,K,F)$ (together with the prime ideal $\wp$)
with parameters defined in Section \ref{section_choice_parameters}. These parameters depend on $n$, $l$, $m$ and  
$\esselle=(s_1,\ldots,s_l) \in \mathcal{L}(\mathtt{r}_1,\ldots,\mathtt{r}_l)$ that we have fixed. Recall that to $\esselle$ we associated a partial
ordering $\prec$ (see Definition \ref{def_prec}) and a relation $\leftrightarrow$ (see Notation 
\ref{notation_lambda_corresponding_lambdal}). Finally, put $s:=s_1+\cdots+s_l$. \fini \\ 

\section{Expression of the matrices $J^{\prec}$ and $J^{\lhd}$} \label{section_expression_J}

\subsection{The matrices $J$ and $J^{\prec}$}

We first give, with our choice of parameters, a simpler expression for $J$. 

\begin{lemma} \label{expression_J} \ \\
\vspace{-7mm} 
\begin{itemize} 
\item[1)] Assume that $(\lambdal,\mul)$ satisfies the conditions $(J_1)$. Then we have
\begin{equation}
j_{\lambdal,\mul} = \left\{ \begin{array}{cl} 
(-1)^{\haut(\rho)+\haut(\rho')} & \mbox{if } \res_n(\hd(\rho)) = \res_n(\hd(\rho')), \\
0 & \mbox{otherwise.}
\end{array} \right.
\end{equation}  
\item[2)] Assume that $(\lambdal,\mul)$ satisfies the conditions $(J_2)$. Then we have
\begin{equation}
j_{\lambdal,\mul} = (-1)^{\haut(\rho)+\haut(\rho')} \varepsilon,
\end{equation}
where
\begin{equation}
\varepsilon := 
\left\{ \begin{array}{cl} 
1 & \mbox{if } \res_n(\hd(\rho)) = \res_n(\hd(\rho')) \quad \mbox{and} \quad \widehat{h} \not \equiv 0 \pmod n, \\
-1 & \mbox{if } \res_n(\hd(\rho)) \neq \res_n(\hd(\rho')) \quad \mbox{and} \quad \widehat{h} \equiv 0 \pmod n, \\
0 & \mbox{otherwise,}
\end{array} \right.
\end{equation}
and $\widehat{h}$ is the common length of $\rho$ and $\rho'$. 
\end{itemize}
\end{lemma}

\begin{demo} 
Let us prove 1). With our choice of parameters, we have 
$$j_{\lambdal,\mul} = (-1)^{\haut(\rho)+\haut(\rho')} \nu_{\wp}(P_{\lambdal,\mul}(x))$$ with
$P_{\lambdal,\mul}(x):=u_d x^{l(j-i)} - u_{d'} x^{l(j'-i')}$. Note that 
$$P_{\lambdal,\mul}(x)=\xi^{a_1}x^{a_2}-\xi^{a_3}x^{a_4},$$ where $\xi \in \mC$ is a primitive $nl$-th root of unity and
$a_1:=dn$, $a_2:=ls_d-dn+l(j-i)$, $a_3:=d'n$ and $a_4:=ls_{d'}-d'n+l(j'-i')$. Using the fact that $\nu_{\wp}(x^N)=0$ for all 
$N \in \mZ$ and $x^N P_{\lambdal,\mul}(x) \in \mC[x]$ for a suitable $N \in \mZ$, we get \vspace{3mm}
$$\begin{array}{l}
\nu_{\wp}(P_{\lambdal,\mul}(x)) \geq 0, \\[3mm] 
\nu_{\wp}(P_{\lambdal,\mul}(x)) \geq 1 \Longleftrightarrow P_{\lambdal,\mul}(\xi)=0  , \\[3mm]
\nu_{\wp}(P_{\lambdal,\mul}(x)) \geq 2 \Longleftrightarrow P_{\lambdal,\mul}(\xi)=P'_{\lambdal,\mul}(\xi)=0.
\end{array}$$

\bigskip \noindent A straightforward computation shows that we have $\nu_{\wp}(P_{\lambdal,\mul}(x)) \geq 1$ if and only if we have
$l(s_d+j-i) \equiv l(s_{d'}+j'-i') \pmod{nl}$, that is if and only if $\res_n(\hd(\rho)) = \res_n(\hd(\rho'))$. Moreover, we have
$$\begin{array}{rcl}
\nu_{\wp}(P_{\lambdal,\mul}(x)) \geq 2 & \Longleftrightarrow & 
 \left\{ \begin{array}{l} a_1+a_2 \equiv a_3+a_4 \pmod{nl} \\ a_2 \xi^{a_1+a_2-1} = a_4 \xi^{a_3+a_4-1} \end{array} \right. \\[4mm]
 & \Longleftrightarrow & \left\{ \begin{array}{l} a_1+a_2 \equiv a_3+a_4 \pmod{nl} \\ a_2 = a_4 \end{array} \right. \\[4mm]
 & \Longleftrightarrow & \left\{ \begin{array}{l} a_1 \equiv a_3 \pmod{nl} \\ a_2 = a_4 \end{array} \right. . \\[4mm] 
\end{array}$$

\noindent But the condition $a_1 \equiv a_3 \pmod{nl}$ implies $d \equiv d' \pmod l$, which is impossible since $d$ and $d'$ are two distinct 
integers ranging from $1$ to $l$. As a consequence, we have $\nu_{\wp}(P_{\lambdal,\mul}(x)) \leq 1$, which proves 1). Let us now prove 2).
With the notation of $(J_2)$, we have $\varepsilon_2=-\varepsilon_1$, whence
$$j_{\lambdal,\mul} = (-1)^{\haut(\rho)+\haut(\rho')} \nu_{\wp}(P_{\lambdal,\mul}(x)), \qquad \mbox{with}$$
$$\begin{array}{rclcl} P_{\lambdal,\mul}(x) &:=& 
\ds \frac{u_d \bigl(x^{l(j-i)} - x^{l(j'-i')} \bigr)}{u_d \bigl(x^{l(j-i)} - x^{l(j''-i'')}\bigr)} \\[5mm] 
&=& \ds \frac{x^{l ((j'-i')-(j-i))}-1}{x^{l ((j''-i'')-(j-i))}-1}
&=& \ds \frac{x^{l \bigl(\cont(\hd(\rho')) - \cont(\hd(\rho)) \bigr)}-1}{x^{l \widehat{h}}-1}\,.
\end{array}$$
In order to complete the proof, we only have to notice that for $N \in \mZ$, we have $\nu_{\wp}(x^{lN}-1)=1$ if $nl$ divides $lN$, that is if
$n$ divides $N$, and $\nu_{\wp}(x^{lN}-1)=0$ otherwise.
\end{demo}

We now start analyzing carefully Cases $(J_1)$ and $(J_2)$. Proposition \ref{analysis_J1_J2_2} gives a characterization in terms of $\lambda$ and 
$\mu$ of the pairs $(\lambdal,\mul)$ that satisfy $(J_1)$ or $(J_2)$. 

\begin{lemma} \label{analysis_J1_J2_1}
Let $\lambdal=(\lambda^{(1)},\ldots,\lambda^{(l)}) \in \Pilm$ and $\mul=(\mu^{(1)},\ldots,\mu^{(l)}) \in \Pilm$ be two distinct multi-partitions, 
and $\lambda$, $\mu \in \Pi$ be such that $\lambdal \leftrightarrow \lambda$ and $\mul \leftrightarrow \mu$.
Then the following statements are equivalent:

 \begin{itemize} 
 \item[\rm (i)] $\lambda \subset \mu$, and $\mu / \lambda$ is a ribbon,
 \item[\rm (ii)] $\lambdal \subset \mul$, and there exists $d \in \interv{1}{l}$ such that $\mu^{(d)} / \lambda^{(d)}$ is a ribbon and 
 $\lambda^{(b)} = \mu^{(b)}$ for all $b \in \interv{1}{l} \setminus \set{d}$.
 \end{itemize}
 \end{lemma}

\begin{demo} 
Let us prove (i) $\Rightarrow$ (ii). By Lemma \ref{ribbon_beta_numbers_1}, passing from $\lambda$ to $\mu$
amounts, as far as abacus diagrams are concerned, to passing from $A(\lambda,s)$ to $A(\mu,s)$ by moving a bead located at position $k$
towards the right. As far as the $l$-runner abacus diagrams $A(\lambdal,\esselle)$ and $A(\mul,\esselle)$ are concerned, 
this amounts to moving a bead located at position $\phi(k)$ on the runner $d:=d(k)$ towards the right. This together with
Lemma \ref{ribbon_beta_numbers_1} applied to $(\lambda^{(d)},\mu^{(d)})$ proves (ii). The converse is similar.   
\end{demo}

Applying twice the previous lemma and Lemma \ref{ribbon_beta_numbers_2} yields the following result.

\begin{prop} \label{analysis_J1_J2_2}
Let $\lambdal$, $\mul \in \Pilm$ be two multi-partitions, and $\lambda$, $\mu \in \Pi$ be such that 
$\lambdal \leftrightarrow \lambda$ and $\mul \leftrightarrow \mu$. Assume that $|\lambda|=|\mu|=r$. Then 
$(\lambdal,\mul)$ satisfies $(J_1)$ or $(J_2)$ if and only if $\sharp \bigl(B(\lambda) \cap B(\mu) \bigr) = r-2$. \cqfd
\end{prop}

The following notation will be very useful.

\begin{notation} \label{notation_H} Let $\lambdal=(\lambda^{(1)},\ldots,\lambda^{(l)}) \in \Pilm$ and $\mul=(\mu^{(1)},\ldots,\mu^{(l)}) \in \Pilm$.
Denote by $\lambda$, $\mu \in \Pi$ the partitions such that $\lambdal \leftrightarrow \lambda$ and $\mul \leftrightarrow \mu$.
Consider the statement
\begin{equation} \label{assumption_H}
\lambdal \prec \mul \quad (\mbox{\ie} \lambda \lhd \mu), \quad |\lambda|=|\mu|=r \quad \mbox{and} 
\quad \sharp \bigl(B(\lambda) \cap B(\mu) \bigr) = r-2.
\end{equation}
Assume now that $(\lambdal, \mul)$ satisfies (\ref{assumption_H}). In this case, we shall use in the sequel the following notation.
Let $(\alpha_1,\ldots,\alpha_r)$ and $(\beta_1,\ldots,\beta_r)$ denote the sequences of $\beta$-numbers associated to $\lambda$ and 
$\mu$ respectively. By Lemma \ref{ribbon_beta_numbers_2} applied to the pair $(\nu,\kappa)=(\lambda,\mu)$, there exist positive integers
$y$, $y'$, $x'$, $x$ and $h$ such that $\set{\alpha_i \mid i \neq x',y'}=\set{\beta_j \mid j \neq x,y}$, 
$\alpha_{y'} = \beta_y-h$ and $\alpha_{x'} = \beta_{x}+h$. Moreover, by Proposition \ref{analysis_J1_J2_2}, $(\lambdal, \mul)$ satisfies $(J_1)$ or 
$(J_2)$. Let $d, d' \in \interv{1}{l}$ and $\widehat{h} \in \mN^{*}$ denote the integers introduced in the definition of Cases $(J_1)$ and 
$(J_2)$ (in Case $(J_2)$, put $d':=d$). Finally, denote by $\gamma$, \resp $\delta \in \interv{0}{nl-1}$ the residue of
$c(\beta_y)-c(\beta_x)$, \resp $n \bigl(d(\beta_y)-d(\beta_x) \bigr)$ modulo $nl$. \fini
\end{notation}

\begin{remark} \label{remark_not_H_implies_j_equals_zero}
Proposition \ref{analysis_J1_J2_2} shows that if $(\lambdal,\mul)$ does not satisfy (\ref{assumption_H}), then 
$j^{\prec}_{\lambdal,\mul}=0$. \fini
\end{remark}

\begin{remark} 
Recall Notation \ref{notation_H} and assume that $(\lambdal,\mul)$ satisfies (\ref{assumption_H}).
Since $\lambda \lhd \mu$, Lemma $\ref{ribbon_beta_numbers_2}$ implies $y' < x'$, so we have  
\begin{equation} \label{eq_1_remark_H_and_beta_numbers}
\beta_x < \beta_x+h = \alpha_{x'} < \alpha_{y'} = \beta_y-h < \beta_y.
\end{equation}
These inequalities, together with $\set{\alpha_i \mid i \neq x',y'}=\set{\beta_j \mid j \neq x,y}$ and the fact that the $\beta_i$'s are
pairwise distinct, imply: 
\begin{equation} \label{eq_2_remark_H_and_beta_numbers}
\set{\beta_x,\beta_y} \cap B(\lambda) = \emptyset.
\end{equation} \fini
\end{remark}

Under assumption (\ref{assumption_H}), the following technical lemma relates some $\beta$-numbers of $\lambda^{(d)}$, $\lambda^{(d')}$,
$\mu^{(d)}$ and $\mu^{(d')}$ on the one hand to some $\beta$-numbers of $\lambda$ and $\mu$ on the other hand.
\bigskip 

\begin{lemma} \label{ribbon_beta_numbers_3}
Recall Notation \ref{notation_H} and assume that $(\lambdal,\mul)$ satisfies (\ref{assumption_H}).

\begin{itemize}
\item[1)] Then, we have 
\begin{equation}
\set{d(\beta_x),d(\beta_y)}=\set{d(\beta_x+h),d(\beta_y-h)}=\set{d,d'}.
\end{equation}

Moreover, for all $1 \leq b \leq l$, we have
\begin{equation} 
\sharp \set{1 \leq i \leq r \mid d(\alpha_i)=b} = \sharp \set{1 \leq i \leq r \mid d(\beta_i)=b}\,;
\end{equation}
let $r_b$ denote this common value. \\
\item[2)] Assume that $(\lambdal, \mul)$ satisfies (\ref{assumption_H}) and $(J_1)$. Then we have $d(\beta_x)=d$ and $d(\beta_y)=d'$. Let
\begin{equation}  
\begin{array}{lcl} 
\boldsymbol{\beta}_{r_d}(\lambda^{(d)})=(\gamma_1,\ldots,\gamma_{r_d}), \quad && \quad
\boldsymbol{\beta}_{r_d}(\mu^{(d)})=(\delta_1,\ldots,\delta_{r_d}), \\
\boldsymbol{\beta}_{r_{d'}}(\lambda^{(d')})=(\gamma'_1,\ldots,\gamma'_{r_{d'}}) \quad & \mbox{and} & \quad
\boldsymbol{\beta}_{r_{d'}}(\mu^{(d')})=(\delta'_1,\ldots,\delta'_{r_{d'}})
\end{array}
\end{equation}
denote the sequences of $\beta$-numbers associated to  $\lambda^{(d)}$, $\mu^{(d)}$, $\lambda^{(d')}$ and $\mu^{(d')}$.
Denote by $b$ (\resp $c$) the row number of the tail (\resp head) of $\rho$, and denote by $b'$ (\resp $c'$) the row number of the tail 
(\resp head) of $\rho'$. By Statement 1), there exist integers $k$, $k'$ such that
\begin{equation} \label{def_k_kprime_ribbon_beta_numbers_3}
\set{k,k'}=\set{\beta_x+h,\beta_y-h}, \quad d(k)=d \quad \mbox{and} \quad d(k')=d'.
\end{equation}
Then we have 
\begin{equation}
\begin{array}{l}
\delta_c=\phi(\beta_x), \quad \gamma_{b}=\phi(k), \quad \gamma'_{c'}=\phi(k'), \quad \delta'_{b'}=\phi(\beta_y) \\[3mm]
\mbox{and} \quad \widehat{h}=\phi(k)-\phi(\beta_x)=\phi(\beta_y)-\phi(k').
\end{array}
\end{equation}
\item[3)] Assume that $(\lambdal, \mul)$ satisfies (\ref{assumption_H}) and $(J_2)$. Denote by
\begin{equation}
\boldsymbol{\beta}_{r_d}(\lambda^{(d)})=(\gamma_1,\ldots,\gamma_{r_d}) \qquad \mbox{and} \qquad
\boldsymbol{\beta}_{r_d}(\mu^{(d)})=(\delta_1,\ldots,\delta_{r_d})
\end{equation}
the sequences of $\beta$-numbers associated to $\lambda^{(d)}$ and $\mu^{(d)}$. Let $b$ (\resp $c$) denote the row number of the tail (\resp head) 
of $\rho$, and $b'$ (\resp $c'$) denote the row number of the tail (\resp head) of $\rho'$. Then we have $\lambda^{(d)} \lhd \mu^{(d)}$ and
\begin{equation}
\begin{array}{l}
\delta_c=\phi(\beta_x), \quad \gamma_b=\phi(\beta_x+h), \quad \gamma_{c'}=\phi(\beta_y-h),\quad \delta_{b'}=\phi(\beta_y) \\[3mm]
\mbox{and} \quad \widehat{h}=\phi(\beta_x+h)-\phi(\beta_x)=\phi(\beta_y)-\phi(\beta_y-h).
\end{array}
\end{equation}
\end{itemize}
\end{lemma}

\begin{demo} 
We pass from $\mul$ to $\lambdal$ by removing the ribbon $\rho'$ and by adding the ribbon $\rho$. This amounts, as far as the abacus diagrams 
$A(\mul,\esselle)$ and $A(\lambdal,\esselle)$ are concerned, to moving two beads  (see the proof of Lemma \ref{analysis_J1_J2_1}). Moving these
two beads amounts, as far as the abacus diagrams $A(\mu,s)$ and $A(\lambda,s)$ are concerned, to moving the beads located at positions
$\set{\beta_x,\beta_y}$ towards the positions $\set{\beta_x+h,\beta_y-h}$, which proves the first two equalities of Statement 1). The last parts 
of Statement 1) come from this and from the equality $\set{\alpha_i \mid i \neq x',y'}=\set{\beta_j \mid j \neq x,y}$. \\
    
\indent Let us now prove Statement 2). A careful analysis of the moves of the beads described above shows more precisely that the following
properties hold:
 
\begin{itemize} \item[(i)]
$\set{\gamma_b,\gamma'_{c'}}=\set{\phi(\beta_x+h),\phi(\beta_y-h)} \qquad \mbox{and} \qquad 
\set{\delta'_{b'},\delta_c}=\set{\phi(\beta_x),\phi(\beta_y)}.$
\item[(ii)] Let $K \in \set{\beta_x,\beta_x+h,\beta_y-h,\beta_y}$. Then we have $\phi(K) \in \set{\gamma_b,\delta_c}$ if and only if 
$d(K)=d$ and $\phi(K) \in \set{\gamma'_{c'},\delta'_{b'}}$ if and only if $d(K)=d'$.
\end{itemize}

Let us first prove that $d(\beta_x)=d$. Assume that $d(\beta_x)=d'$. By (i) and (ii), we have 
$\phi(\beta_x) \in \set{\delta'_{b'},\delta_c} \cap \set{\gamma'_{c'},\delta'_{b'}}$. This implies $\phi(\beta_x)=\delta'_{b'}$. Indeed, if
$\phi(\beta_x) \neq \delta'_{b'}$, we must have $\phi(\beta_x)=\delta_c=\gamma'_{c'} \in \set{\phi(\beta_x+h),\phi(\beta_y-h)}$ by (i).
Let $K \in \set{\beta_x+h,\beta_y-h}$ be such that $\phi(K)=\gamma'_{c'}=\phi(\beta_x)$. By (ii), we have $d(K)=d(\beta_x)$; moreover, 
we have $\phi(K)=\phi(\beta_x)$, whence $K=\beta_x$. This contradicts (\ref{eq_1_remark_H_and_beta_numbers}),    
so $\phi(\beta_x)=\delta'_{b'}$. Let $K \in \set{\beta_x+h,\beta_y-h}$ be such that $d(K)=d'$. By (\ref{eq_1_remark_H_and_beta_numbers}) and 
(\ref{property_P2}), we have $\phi(K)>\phi(\beta_x)=\delta'_{b'}$, hence by (ii) we have $\phi(K)=\gamma'_{c'}$. As a 
consequence, we have $\gamma'_{c'} > \delta'_{b'}$. Moreover, by Lemma \ref{ribbon_beta_numbers_1} applied to $(\lambda^{(d')},\mu^{(d')})$, we have 
$\gamma'_{c'}=\delta'_{b'}-\widehat{h} < \delta'_{b'}$, which is absurd. By Statement 1), we thus have $d(\beta_x)=d$ and $d(\beta_y)=d'$. 
Now let $k$ be the integer defined by (\ref{def_k_kprime_ribbon_beta_numbers_3}) and assume that $k=\beta_x+h$ (the proof for the case $k=\beta_y-h$
is similar). By (ii) we have $\set{\gamma_b,\delta_c}=\set{\phi(\beta_x),\phi(\beta_x+h)}$. Moreover, by Lemma \ref{ribbon_beta_numbers_1} applied 
to $(\lambda^{(d)},\mu^{(d)})$, we have $\delta_c=\gamma_b-\widehat{h} < \gamma_b$. By (\ref{eq_1_remark_H_and_beta_numbers}) and 
(\ref{property_P2}), we therefore have 
$\delta_c = \phi(\beta_x)$ and $\gamma_b = \phi(\beta_x+h)$, whence $\widehat{h}= \gamma_b - \delta_c = \phi(\beta_x+h) - \phi(\beta_x)$. By a 
similar argument, we get $\gamma'_{c'}=\phi(\beta_y-h)$, $\delta'_{b'}=\phi(\beta_y)$ and $\widehat{h}= \phi(\beta_y) - \phi(\beta_y-h)$. \\

\indent Let us now prove 3). Since $\lambda \lhd \mu$, by (\ref{eq_1_remark_H_and_beta_numbers}) and (\ref{property_P2}) we have 
$$\phi(\beta_x) < \phi(\beta_x+h) < \phi(\beta_y-h) < \phi(\beta_y) .$$ Moreover, a careful analysis of the 
moves of the beads mentioned at the beginning of the proof shows that
$$\set{\delta_{b'},\delta_c}=\set{\phi(\beta_x),\phi(\beta_y)} \qquad \mbox{and} 
\qquad \set{\gamma_{b},\gamma_{c'}}=\set{\phi(\beta_x+h),\phi(\beta_y-h)}.$$ Assume that $\delta_c = \phi(\beta_y)$. Since
$\phi(\beta_y) > \phi(\beta_x+h)$, $\phi(\beta_y) > \phi(\beta_y-h)$ and $\gamma_{b}$ is in the set $\set{\phi(\beta_x+h),\phi(\beta_y-h)}$, 
we must have $\delta_c > \gamma_{b}$. Moreover, applying Lemma \ref{ribbon_beta_numbers_2} to the pair $(\lambda^{(d)},\mu^{(d)})$ yields
$\gamma_b=\delta_c+\widehat{h} > \delta_c$, which is absurd. We thus have $\delta_c=\phi(\beta_x)$ and $\delta_{b'}=\phi(\beta_y)$. Since 
$\phi(\beta_x)<\phi(\beta_y)$, we have $\delta_c<\delta_{b'}$, whence $c>b'$. Applying again Lemma \ref{ribbon_beta_numbers_2} shows that 
$\lambda^{(d)} \lhd \mu^{(d)}$ and $b' \leq c' < b \leq c$. In particular, we have $c'<b$, whence $\gamma_b<\gamma_{c'}$. Since  
$\phi(\beta_x+h)<\phi(\beta_y-h)$, we have $\gamma_b=\phi(\beta_x+h)$ and $\gamma_{c'} = \phi(\beta_y-h)$. Lemma \ref{ribbon_beta_numbers_2} then 
implies that $\widehat{h} = \gamma_{b} - \delta_c = \phi(\beta_x+h)-\phi(\beta_x)$ and 
$\widehat{h} = \delta_{b'} - \gamma_{c'} = \phi(\beta_y)-\phi(\beta_y-h)$.  
\end{demo}

We are now ready to derive the expressions of the $j^{\prec}_{\lambdal,\mul}$'s that we need for proving Theorem \ref{thm_A_J}.

\begin{prop} \label{proposition_J1_J2} Recall Notation \ref{notation_H}. 

\begin{itemize}
\item[1)] Assume that $(\lambdal,\mul)$ satisfies $(J_1)$ and $\lambdal \prec \mul$. Then we have $\delta>0$ and 
\begin{equation}
j^{\prec}_{\lambdal,\mul}=
\left\{ \begin{array}{cl} 
(-1)^{\haut(\rho)+\haut(\rho')} & \mbox{if } h \equiv \gamma \pmod{nl} \quad \mbox{or} \quad h \equiv \delta \pmod{nl}, \\
0 & \mbox{otherwise.}
\end{array} \right.
\end{equation}

\item[2)] Assume that $(\lambdal,\mul)$ satisfies $(J_2)$ and $\lambdal \prec \mul$. Then we have $\delta=0$ and 
\begin{equation}
j^{\prec}_{\lambdal,\mul} = (-1)^{\haut(\rho)+\haut(\rho')} \varepsilon,
\end{equation}
where 
\begin{equation}
\varepsilon := \left\{ \begin{array}{cl} 
1 & \mbox{if } h \equiv \gamma \pmod{nl} \quad \mbox{and} \quad h \not \equiv 0 \pmod{nl}, \\
-1 & \mbox{if } h \not \equiv \gamma \pmod{nl} \quad \mbox{and} \quad h \equiv 0 \pmod{nl}, \\
0 & \mbox{otherwise.}
\end{array} \right.
\end{equation}
\end{itemize}
\end{prop}

\pv 
We prove only Statement 1), the proof of Statement 2) being similar. Recall the notation from Lemma \ref{ribbon_beta_numbers_3}. The statement 
$\delta>0$ comes from Statement 1) of that lemma. Applying Lemma \ref{ribbon_beta_numbers_1} to the pairs $(\lambda^{(d)},\mu^{(d)})$ and 
$(\lambda^{(d')},\mu^{(d')})$ yields 
$$\res_n(\hd(\rho)) = \delta_c \bmod n \qquad \mbox{and} \qquad \res_n(\hd(\rho')) = \gamma'_{c'} \bmod n .$$ 
By Lemma \ref{expression_J}, it is thus enough to prove the equivalence
$$\delta_c \equiv \gamma'_{c'} \pmod n \Longleftrightarrow 
\bigl(h \equiv \gamma \pmod{nl} \quad \mbox{ or } \quad h \equiv \delta \pmod{nl} \bigr) .$$ 
By Statement 1) of Lemma \ref{ribbon_beta_numbers_3}, one of the two following cases occurs. 

\begin{itemize}
\item \underline{\emph{First case:} we have $d(\beta_x)=d(\beta_x+h)$ and $d(\beta_y)=d(\beta_y-h)$}. By Lemma \ref{ribbon_beta_numbers_3}
and (\ref{property_P1}), we have the following equivalences, where congruences stand modulo $n$: 
$$\delta_c \equiv \gamma'_{c'} \Longleftrightarrow \phi(\beta_x) \equiv \phi(\beta_y-h) 
\Longleftrightarrow \beta_x \equiv \beta_y-h \Longleftrightarrow h \equiv \gamma .$$
It remains thus to prove that $\delta_c \equiv \gamma'_{c'} \pmod n \Longrightarrow h \equiv \gamma \pmod{nl}$. Assume that 
$\delta_c \equiv \gamma'_{c'} \pmod n$; we then have $h \equiv \gamma \pmod{n}$. This and the equality $d(\beta_x)=d(\beta_x+h)$ force 
$d(h)=1$, whence $h \equiv \gamma \pmod{nl}$. \\

\item \underline{\emph{Second case:} we have $d(\beta_x)=d(\beta_y-h)$ and $d(\beta_y)=d(\beta_x+h)$}. By arguing as above we prove the 
equivalence $$\delta_c \equiv \gamma'_{c'} \pmod n \Longleftrightarrow h \equiv 0 \pmod n .$$
Assume that $\delta_c \equiv \gamma'_{c'} \pmod{n}$. Then we have $h \equiv 0 \pmod n$, whence 
$$d(\beta_y)=d(\beta_x+h) \equiv d(\beta_x) + d(h) -1 \pmod l,$$ whence 
$n(d(h)-1) \equiv n \bigl(d(\beta_y)-d(\beta_x) \bigr) \equiv \delta \pmod{nl}$ and $h \equiv \delta \pmod{nl}$. 
As a consequence, we have in this case $\delta_c \equiv \gamma'_{c'} \pmod{n} \Longleftrightarrow h \equiv \delta \pmod{nl}$, 
which completes the proof. \cqfd
\end{itemize} 

\begin{remark} \label{remark_proof_proposition_J1_J2}
Recall Notation \ref{notation_H} and assume that $(\lambdal,\mul)$ satisfies $(J_1)$ and $j^{\prec}_{\lambdal,\mul} \neq 0$.
Then the proof of Proposition \ref{proposition_J1_J2} shows in particular that 
\begin{equation}
\phi(\beta_x) \equiv \phi(k') \pmod n,
\end{equation}
where $k' \in \mZ$ is defined by (\ref{def_k_kprime_ribbon_beta_numbers_3}). \fini
\end{remark}

\subsection{What happens if the multi-charge $\esselle$ is $m$-dominant} \label{section_case_sl_dominant}

The goal of this section is to show that with our choice of parameters, the matrix $J^{\lhd}$ is a special case of a matrix $J^{\prec}$ 
when the multi-charge $\esselle$ is $m$-dominant (see Proposition \ref{expr_J_dominant}). However, the results we prove here will not be used for 
the proof of Theorem \ref{thm_A_J}. 

\begin{lemma} \label{lemma_weight_partitions}
Let $\lambdal,  \mul \in \Pilm$, and let $\lambda$, $\mu \in \Pi$ be such that $\lambdal \leftrightarrow \lambda$ and $\mul \leftrightarrow \mu$. 
Assume that $j_{\lambdal,\mul} \neq 0$. Then we have, with the notation from $(J_1)$, $(J_2)$ and (\ref{notation_nb_i_nodes})\,: 

\begin{equation}
N_i(\rho)=N_i(\rho') \quad (i \in \mZ) \qquad \mbox{and} \qquad |\lambda|=|\mu|.
\end{equation} 
\end{lemma} 

\begin{demo} 
Since $j_{\lambdal,\mul}$ is nonzero, Lemma \ref{expression_J} shows that at least one of the following cases occurs:

\begin{itemize}
\item \underline{\emph{First case:}} we have $\res_n(\hd(\rho))=\res_n(\hd(\rho'))$. Note that if $\rho$ is a ribbon, then the integers  
$\cont(\gamma),  \gamma \in \rho$ are pairwise distinct, and the set formed by these numbers is exactly the interval 
$\interv{\cont({\hd(\rho)})}{\cont{(\tl(\rho)})}$. Combining this with the assumption $\res_n(\hd(\rho))=\res_n(\hd(\rho'))$
and $\widehat{h}=\ell(\rho)=\ell(\rho')$, we get that $N_i(\rho)=N_i(\rho')$ for any $i \in \mZ$.
\item \underline{\emph{Second case:}} we have $\widehat{h} \equiv 0 \pmod n$. Then for any $i \in \mZ$, we have the equalities
$N_i(\rho)=N_i(\rho')=\widehat{h}/n$.
\end{itemize}

Let us now show that $|\lambda|=|\mu|$. Let $\boldsymbol{\nu}_l:=\lambdal \cap \mul$, and let $\nu \in \Pi$ be such that 
$\boldsymbol{\nu}_l \leftrightarrow \nu$. We claim that
$$h:=|\lambda|-|\nu|=\bigl( (n-1)l + 1 \bigr) N_0(\rho) + \bigl(\widehat{h}-N_0(\rho) \bigr).$$
By induction on $\widehat{h}$, we can restrict ourselves to the case when $\widehat{h}=1$, that is $\rho$ contains a single node $\gamma$. 
Let $r \in \mN$ be such that $\lambda$ and $\nu$ have at most $r$ parts. By Lemma \ref{ribbon_beta_numbers_1}, there exist 
$\alpha \in B_r(\nu)$ and $\beta \in B_r(\lambda)$ such that 
$B_r(\nu) \setminus \set{\alpha}=B_r(\lambda) \setminus \set{\beta}$ and
$\alpha = \beta - h$. The abacus diagrams $A(\nu,s)$ and $A(\lambda,s)$ differ only by the moving of a bead; the same thing holds for
the diagrams $A(\boldsymbol{\nu}_l,\esselle)$ and $A(\lambdal,\esselle)$. By considering the initial and the final positions of these two beads, 
we get $$\phi(\beta)=\phi(\alpha)+\widehat{h}=\phi(\alpha)+1 \qquad \mbox{and} \qquad d(\beta)=d(\alpha)=d.$$      
Moreover, by Lemma \ref{ribbon_beta_numbers_1} and (\ref{property_P1}), we have $\res_n(\gamma)=\phi(\alpha) \bmod n = \alpha \bmod n$. 
Let us now distinguish two cases. If $\res_n(\gamma)=0 \bmod n$ (\ie if $N_0(\rho)=1$), then we have $\alpha = n + n(d-1) + nlm$ with $m \in \mZ$, 
whence $\phi(\beta)=\phi(\alpha)+1=1+n(m+1)$. Since $d(\beta)=d$, we get $\beta=1+n(d-1)+nl(m+1)$, whence $h=\beta-\alpha=(n-1)l+1$. 
Similarly, if $\res_n(\gamma) \neq 0 \bmod n$ (\ie if $N_0(\rho)=0$), then we have $h=1$. This proves the claimed formula.
In a similar way we prove that $|\mu|-|\nu|=\bigl( (n-1)l + 1 \bigr) N_0(\rho') + \bigl(\widehat{h}-N_0(\rho') \bigr)$. Since
$N_0(\rho)=N_0(\rho')$, we do have $|\lambda|=|\mu|$.
\end{demo}

\begin{lemma} \label{lemma_coincidence_orderings}
Let $\lambdal,  \mul \in \Pilm$, and $\lambda$, $\mu \in \Pi$ be such that $\lambdal \leftrightarrow \lambda$ and 
$\mul \leftrightarrow \mu$. Assume that $|\lambda|=|\mu|$. Consider the following cases:

\begin{itemize}
\item[1)] $(\lambdal,\mul)$ satisfies $(J_1)$, $j_{\lambdal,\mul} \neq 0$ and $\esselle$ is $m$-dominant,
\item[2)] $(\lambdal,\mul)$ satisfies $(J_2)$.
\end{itemize}

Then in either case, we have: $\lambdal \prec \mul \Longleftrightarrow \lambdal \lhd \mul$.
\end{lemma}

\begin{demo} 
In either case, we can apply Proposition \ref{analysis_J1_J2_2} and then Lemma \ref{ribbon_beta_numbers_2} to get that 
$\lambdal \prec \mul$ or $\mul \prec \lambdal$. It is thus enough to prove that $\lambdal \prec \mul \Rightarrow \lambdal \lhd \mul$. 
Assume from now on that $\lambdal \prec \mul$ and $(\lambdal,\mul)$ satisfies either case of Lemma \ref{lemma_coincidence_orderings}. Then by 
Proposition \ref{analysis_J1_J2_2}, $(\lambdal,\mul)$ satisfies (\ref{assumption_H})\,; therefore (\ref{eq_1_remark_H_and_beta_numbers}) holds. 
Recall Notation \ref{notation_H}. If $(\lambdal,\mul)$ satisfies $(J_2)$, then by Statement 3) of Lemma \ref{ribbon_beta_numbers_3}, we have 
$\lambda^{(d)} \lhd \mu^{(d)}$. Moreover, for all $b \in \interv{1}{l} \setminus \set{d}$ we have $\lambda^{(b)} = \mu^{(b)}$, whence 
$\lambdal \lhd \mul$. Assume now that $(\lambdal,\mul)$ satisfies $(J_1)$, $j_{\lambdal,\mul} \neq 0$ and $\esselle$ is $m$-dominant. The key point 
of the proof is the following. Let $\boldsymbol{\nu}_l=\bigl( \nu^{(1)},\ldots,\nu^{(l)} \bigr) \in \Pilm$
be such that $|\nu|=r$, where $\nu$ is the partition such that $\boldsymbol{\nu}_l \leftrightarrow \nu$. Then under the assumption that 
$\esselle$ is $m$-dominant, we have 
\begin{equation} \bigl( d(k)<d(k'), \ k,  k' \in B(\nu) \bigr) \Rightarrow \phi(k) \geq \phi(k').
\tag{$*$}  
\end{equation}
Indeed, let $k, k' \in B(\nu)$, $b:=d(k)$, $b':=d(k')$ and $N$ (\resp $N'$) be the number of parts of $\nu^{(b)}$ (\resp $\nu^{(b')}$). Since 
$\boldsymbol{\nu}_l \leftrightarrow \nu$, we have $\phi(k) \in B_N(\nu^{(b)})$ and $\phi(k') \in B_{N'}(\nu^{(b')})$. As a consequence, there exist 
$i \in \interv{1}{N}$, $i' \in \interv{1}{N'}$ such that $\phi(k)=s_b+\nu^{(b)}_i-i+1$ and $\phi(k')=s_{b'}+\nu^{(b')}_{i'}-i'+1$. Since 
$\esselle$ is $m$-dominant and $b<b'$, we have 
$$\begin{array}{rl}
\phi(k)-\phi(k') & = (s_b-s_{b'}) + (i'-i) + (\nu^{(b)}_i-\nu^{(b')}_{i'}) \\[3mm]
& \geq s_b-s_{b'}-N-\nu^{(b')}_{i'} \quad\geq\quad s_b-s_{b'} - (|\nu^{(b)}|+|\nu^{(b')}|) \quad\geq\quad s_b-s_{b'}-|\boldsymbol{\nu}_l| \\[3mm]
&\geq 0,
\end{array}$$  
which shows $(*)$. We now claim that $d(\beta_x) \geq d(\beta_y)$. Assume indeed that $d(\beta_x) < d(\beta_y)$. Recall that 
(\ref{eq_1_remark_H_and_beta_numbers}) holds. By $(*)$ applied to $(k,k')=(\beta_x,\beta_y)$, we have $\phi(\beta_x) \geq \phi(\beta_y)$. 
This, (\ref{eq_1_remark_H_and_beta_numbers}) and (\ref{property_P3}) imply that $m(\beta_x)=m(\beta_y)$, where the map $k \mapsto m(k)$ is
defined in Section \ref{section_notation_tau_l}. Since this map is increasing, we have by (\ref{eq_1_remark_H_and_beta_numbers})\,: 
$$m(\beta_x) \leq m(\beta_x+h) \leq m(\beta_y-h) \leq m(\beta_y)=m(\beta_x),$$ so 
equalities hold throughout. Let now $k' \in \mZ$ be the integer defined by (\ref{def_k_kprime_ribbon_beta_numbers_3}). Since 
$j^{\prec}_{\lambdal,\mul} \neq 0$ by assumption, we can apply Remark \ref{remark_proof_proposition_J1_J2}  
and get $\phi(\beta_x) \equiv \phi(k') \pmod n$. This together with $m(\beta_x)=m(k')$ forces $\phi(\beta_x)=\phi(k')$. Moreover, by Statement 2 of
Lemma \ref{ribbon_beta_numbers_3}, we have $0 < \widehat{h}=\phi(\beta_y)-\phi(k') \leq \phi(\beta_x)-\phi(k')$, which is absurd. Therefore we have 
$d(\beta_x) \geq d(\beta_y)$ as claimed. Again by Statement 2) of Lemma \ref{ribbon_beta_numbers_3}, we have $d(\beta_x)=d$ and $d(\beta_y)=d'$, 
whence $d>d'$. Moreover, we have $|\mu^{(d')}| = |\lambda^{(d')}|+\widehat{h}$, $|\mu^{(d)}| = |\lambda^{(d)}|-\widehat{h}$ and 
$\mu^{(b)} = \lambda^{(b)}$ for all $b \in \interv{1}{l} \setminus \set{d,d'}$. This and the inequality $d'<d$ imply  
$\lambdal \lhd \mul$.
\end{demo}

\begin{prop} \label{expr_J_dominant}
Assume that $\esselle$ is $m$-dominant. Then we have $J^{\prec}=J^{\lhd}$.
\end{prop} 

\begin{demo}
Let $\lambdal$, $\mul \in \Pilm$. If $j_{\lambdal,\mul}=0$, then we have $j^{\prec}_{\lambdal,\mul} = j^{\lhd}_{\lambdal,\mul}=0$ and we are done.
Assume now that $j_{\lambdal,\mul} \neq 0$. It is enough to prove that $\lambdal \prec \mul \Leftrightarrow \lambdal \lhd \mul$. Note that by Lemma 
\ref{lemma_weight_partitions}, we have $|\lambda|=|\mu|$, where $\lambda$, $\mu \in \Pi$ are such that $\lambdal \leftrightarrow \lambda$ and 
$\mul \leftrightarrow \mu$. Moreover, since $j_{\lambdal,\mul} \neq 0$, $(\lambdal,\mul)$ satisfies either $(J_1)$ or $(J_2)$. We can therefore 
apply Lemma \ref{lemma_coincidence_orderings} to conclude. 
\end{demo}

\begin{remark} 
The reader should be warned that the orderings $\prec$ and $\lhd$ do not necessarily coincide, even if the multi-charge $\esselle$ is $m$-dominant. 
For example, let $n = 2$, $l = 2$, $m=6$, $\esselle = (3,-3)$, $\lambdal=\bigl( (2,1),(1,1,1) \bigr)$
and $\mul=\bigl( (3),(2,1) \bigr)$. Then we have $\lambdal \lhd \mul$; however, the partitions $\lambda$
and $\mu$ such that $\lambdal \leftrightarrow \lambda$ and $\mul \leftrightarrow \mu$ are $\lambda=(9,6,3,1,1,1,1,1,1,1)$ and 
$\mu=(10,3,3,2,2,2,1,1,1)$, so $\lambdal$ and $\mul$ are not comparable with respect to $\preceq \;$. \fini
\end{remark}

\section{Admissible sequences, good sequences} \label{section_good_sequence}

In this section we compute the matrix $A'(1)$. To this aim, we examine in detail the straightening of the wedge product 
$v_{\textbf{k}} = v_{k_1} \wedge \cdots \wedge v_{k_r}$. If $v_{\textbf{k}}$ is not 
ordered, there are in general several ways to straighten it by applying recursively the rules $(R_1)$-$(R_4)$. In the sequel, we decide to 
straighten at each step the \emph{first} infraction that occurs in $v_{\textbf{k}}$, that is, the first 
$v_{k_i} \wedge v_{k_{i+1}}$ with $k_i \leq k_{i+1}$. This leads to the notion of admissible sequence that we introduce in Definition
\ref{def_admissible_sequence}. Fix an entry $a_{\lambdal,\mul}'(1)$ of $A'(1)$. We give in Proposition \ref{express_aprime_lambda_mu} an 
expression of it in terms of admissible sequences. Each sequence having a nonzero contribution is called a good sequence. We then show that
there exists at most one good sequence (see Propositions \ref{prop_existence_good_sequence} and \ref{prop_uniqueness_good_sequence}), and it if
exists we compute its length modulo $2$ (see Proposition \ref{compatibility_signs_thm_A_J}).

\subsection{Definitions} 

\begin{definition} \label{def_adjacent_wedges} 
Let $\textbf{k}=(k_1,\ldots,k_r)$ and $\textbf{l}=(l_1,\ldots,l_r) \in \mZ^{r}$. We say that the wedge products 
$v_{\textbf{k}}$ and $v_{\textbf{l}}$ are \emph{adjacent} if there exists $1 \leq i \leq r-1$ ($i$ is then necessarily unique) such that: \\ 
\begin{itemize}
\item[(i)] $k_i \leq k_{i+1}$, and $ k_j >k_{j+1}$ for all $1 \leq j \leq i-1$, 
\item[(ii)] $k_j = l_j$ for all $j \in \interv{1}{r} \setminus \set{i,i+1}$,
\item[(iii)] the wedge product $v_{l_i} \wedge v_{l_{i+1}}$ appears in the straightening of $v_{k_i} \wedge v_{k_{i+1}}$.
\end{itemize} 

In this case, denote by $t \in \interv{1}{4}$ the index of the rule $(R_t)$ applied for the straightening of $v_{k_i} \wedge v_{k_{i+1}}$ and
by $\alpha(v_{\textbf{k}},v_{\textbf{l}}) \in \mZ[q,q^{-1}]$ the coefficient of $v_{l_i} \wedge v_{l_{i+1}}$ in the resulting linear combination.
If $(l_i,l_{i+1})=(k_{i+1},k_i)$, then write $v_{\textbf{k}} \stackrel{\bullet}{\rightarrow} v_{\textbf{l}}$ and set
$m(v_{\textbf{k}},v_{\textbf{l}}):=0$. Otherwise, write $v_{\textbf{k}} \stackrel{t}{\rightarrow} v_{\textbf{l}}$ and set 
$m(v_{\textbf{k}},v_{\textbf{l}}):=1$; note that $t \geq 2$ in this case. In either case, write more simply 
$v_{\textbf{k}} \rightarrow v_{\textbf{l}}$. 
\end{definition}
 
\begin{definition} \label{def_admissible_sequence}

The sequence $\textbf{V} = (v_{\textbf{k}_i})_{0 \leq i \leq N}$ is called \emph{admissible} if each $v_{\textbf{k}_i}$ is a wedge product of
$r$ factors and if we have 
\begin{equation}
v_{\textbf{k}_0} \rightarrow v_{\textbf{k}_1} \rightarrow \cdots \rightarrow v_{\textbf{k}_N} \, ;
\end{equation}
in this case, $N$ is called the \emph{length} of the sequence $\textbf{V}$. Set 
\begin{equation}
\alpha_{\textbf{V}}(q):=\prod_{i = 1}^{N}{\alpha(v_{\textbf{k}_{i-1}},v_{\textbf{k}_i})} \in \mZ[q, q^{-1}] \qquad \mbox{and} \qquad
m(\textbf{V}):= \sum_{i=1}^{N}{m(v_{\textbf{k}_{i-1}},v_{\textbf{k}_i})} \in \mN.
\end{equation}

\medskip

Recall Notation \ref{notation_permuted_wedges} and the definition of $\omega$ from Notation \ref{notation_S_r}. Let $\lambda$, $\mu \in \Pi$ be two 
partitions of $r$. We say that the sequence of wedge products $\textbf{V}=(v_{\textbf{k}_i})_{0 \leq i \leq N}$ is \emph{$(\lambda,\mu)$-admissible} 
if it is an admissible sequence of wedge products (of $r$ factors for each of them) such that $v_{\textbf{k}_N}=v_{\lambda}$ and 
$v_{\textbf{k}_0} = v_{\omega.\mu}$.~\fini
\end{definition}

\begin{remark} It is easy to see that if $\lambda \neq \mu$, then there cannot exist any $(\lambda,\mu)$-admissible sequence 
$\textbf{V}=(v_{\textbf{k}_i})_{0 \leq i \leq N}$ such that $m(\textbf{V})=0$. \fini
\end{remark}

\begin{definition} \label{def_good_sequence}
A $(\lambda,\mu)$-admissible sequence $\textbf{V}$ such that $m(\textbf{V})=1$ is called a \emph{good} sequence (with respect to 
$(\lambda,\mu)$). Such a sequence can be written as  
\begin{equation}
v_{\omega.\mu} \stackrel{\bullet}{\rightarrow} \cdots \stackrel{\bullet}{\rightarrow} \textbf{u} 
\stackrel{t}{\rightarrow} \textbf{v} \stackrel{\bullet}{\rightarrow} \cdots \stackrel{\bullet}{\rightarrow} v_{\lambda},
\end{equation} 
with $t \in \interv{2}{4}$. \fini
\end{definition}

\begin{remark} \label{rq_CN_existence_bonne_sequence}
If a good sequence (with respect to $(\lambda,\mu)$) exists, then we have 
\begin{equation}
\sharp \bigl(B(\lambda) \cap B(\mu) \bigr) = r-2.
\end{equation}
\fini
\end{remark}

\subsection{Reduction to the good sequences}

Recall the expression of the involution $\barre$ of $\Fq_m$ given in (\ref{eq_action_barre}).
Expressing this involution in terms of the $v_{\textbf{k}}$'s and then using Definition 
\ref{def_admissible_sequence} yields the following expression for the coefficients of $A(q)$. 

\begin{lemma} \label{matrix_A_admissible_sequences}
Let $\lambdal, \ \mul \in \Pilm$ and $\lambda$, $\mu \in \Pi$ be the partitions such that $\lambdal \leftrightarrow \lambda$ and
$\mul \leftrightarrow \mu$. Assume that $|\lambda|=|\mu|$. Then we have \em
\begin{equation}
a_{\lambdal,\mul}(q)= \varepsilon(\lambdal,\mul) \, q^{\kappa(\textbf{d}(\mu))-\kappa(\textbf{c}(\mu))}
\sum_{\textbf{V}}{\alpha_{\textbf{V}}(q)},
\end{equation}
\em where the sum ranges over all $(\lambda,\mu)$-admissible sequences $\emph{\textbf{V}}$, and 
$\varepsilon(\lambdal,\mul)$ is the sign defined~by \em
\begin{equation} \label{eq_sign_epsilon_lambdal_mul}
\varepsilon(\lambdal,\mul):=(-1)^{ \kappa(\textbf{d}(\mu))+\ell(v(\omega. \mu))+\ell(v(\lambda))}.
\end{equation}
\cqfd
\end{lemma}

\begin{remark} One should be aware that in general several terms might contribute to this sum, so the statement at the beginning of 
\cite[Section 4]{Ry} is not correct. However, we can fix the argument from \cite{Ry} by showing first that only good sequences do contribute to
$a'_{\lambdal,\mul}(1)$ (see Proposition \ref{express_aprime_lambda_mu}), and then that there exists at most one good sequence (see Propositions 
\ref{prop_existence_good_sequence} and \ref{prop_uniqueness_good_sequence}). \fini
\end{remark}

\begin{prop} \label{express_aprime_lambda_mu} Let $\lambdal, \ \mul \in \Pilm$ be two distinct multi-partitions, and
$\lambda$, $\mu \in \Pi$ be such that $\lambdal \leftrightarrow \lambda$ and $\mul \leftrightarrow \mu$. Assume that
$|\lambda|=|\mu|=r$. Then we have \em
\begin{equation}
a'_{\lambdal,\mul}(1)= \varepsilon(\lambdal,\mul) \sum_{\textbf{V}}{\alpha'_{\textbf{V}}(1)},
\end{equation} 
\em where the sum ranges over all the good sequences with respect to $(\lambda,\mu)$, and $\varepsilon(\lambdal,\mul)$ is the sign defined by 
(\ref{eq_sign_epsilon_lambdal_mul}).
\end{prop}

\begin{demo} 
By Lemma \ref{matrix_A_admissible_sequences}, we have $$a_{\lambdal,\mul}(q)=\varepsilon(\lambdal,\mul) \sum_{\textbf{V}}{f_{\textbf{V}}(q)} ,$$
where the sum ranges over all the $(\lambda,\mu)$-admissible sequences $\textbf{V}$ and $f_{\textbf{V}}(q)$ is the Laurent polynomial 
defined by $f_{\textbf{V}}(q) := q^{\kappa(\textbf{d}(\mu))-\kappa(\textbf{c}(\mu))} \alpha_{\textbf{V}}(q)$. Note that if 
$\textbf{V}$ is an admissible sequence, then we have $m(\textbf{V}) \geq 1$ (because $\lambda \neq \mu$), and moreover the rules $(R_1)$-$(R_4)$ 
imply that $$\alpha_{\textbf{V}}(q) \in (q^2-1)^{m(\textbf{V})} \mZ[q,q^{-1}] .$$ As a consequence, if $\textbf{V}$ is a $(\lambda,\mu)$-admissible
sequence such that $m(\textbf{V}) \geq 2$, then $(q^2-1)^2$ divides $f_{\textbf{V}}(q)$ in $\mZ[q,q^{-1}]$, whence $f_{\textbf{V}}'(1)=0$.
Moreover, if $\textbf{V}$ is a good sequence, then the previous discussion shows that $ \alpha_{\textbf{V}}(1)=0$, whence
$$f_{\textbf{V}}'(1)=\bigl(\kappa(\textbf{d}(\mu))-\kappa(\textbf{c}(\mu)) \bigr) \alpha_{\textbf{V}}(1) + \alpha_{\textbf{V}}'(1) 
= \alpha_{\textbf{V}}'(1).$$   
\end{demo}

\subsection{Existence and uniqueness of the good sequence}

We first give (see Proposition \ref{prop_existence_good_sequence}) some sufficient conditions for the existence of a good sequence 
(with respect to a given pair $(\lambda,\mu)$). In order to do this, we must study the sequence of permutations that we apply to the components 
of wedge products when we go through an admissible sequence
\begin{equation}
v_{\textbf{k}_0} \stackrel{\bullet}{\rightarrow} v_{\textbf{k}_1} \stackrel{\bullet}{\rightarrow} \cdots \stackrel{\bullet}{\rightarrow} 
v_{\textbf{k}_N},
\end{equation}
where $k_0,\ldots,k_N \in \mZ^r$ and $v_{\omega.\textbf{k}_0}$ is ordered. 
\begin{notation}
Let $\sigma$, $\tau \in \Sr$. Write 
\begin{equation}
\sigma \rightarrow \tau
\end{equation}
if there exists $1 \leq i \leq r-1$ such that $\sigma(i)<\sigma(i+1)$, 
\emph{with $i$ minimal for this property}, and such that $\tau = \sigma_i \sigma$. The relation $\rightarrow$ on $\Sr$ is closely related to the
relation $\stackrel{\bullet}{\rightarrow}$ on wedge products of $r$ factors. Recall that for $v_{\textbf{k}}=v_{k_1} \wedge \cdots \wedge v_{k_r}$ 
and $\sigma \in \Sr$, we have $v_{\sigma.\textbf{k}}=v_{k_{\sigma^{-1}(1)}} \wedge  \cdots \wedge v_{k_{\sigma^{-1}(r)}}$. Then by definition we 
have $\sigma \rightarrow \tau$ if and only if $v_{\sigma^{-1}.\textbf{k}} \stackrel{\bullet}{\rightarrow} v_{\tau^{-1}.\textbf{k}}$, where 
$\textbf{k}=(k_1,\ldots,k_r) \in \mZ^r$ is such that $v_{\omega.\textbf{k}}$ is ordered (namely, $k_1<\cdots<k_r$). \fini
\end{notation}

\bigskip

Consider now the following reduced expression for the longest element in $\Sr$ : 
\begin{equation}
\omega = (\sigma_1\sigma_2 \cdots \sigma_{r-1})(\sigma_1\sigma_2 \cdots \sigma_{r-2}) \cdots(\sigma_1\sigma_2)(\sigma_1), 
\end{equation}

\noindent and for $0 \leq i \leq \frac{r(r-1)}{2}$ let $\omega[i]$ denote the right factor of length $i$ in this word (by convention,
$\omega[0]=\mathrm{id}$). For example, for $r \geq 3$ we have $\omega[5]=\sigma_2\sigma_3\sigma_1\sigma_2\sigma_1$. The sequence 
$(\omega[i])_{0 \leq i \leq \frac{r(r-1)}{2}}$ enjoys the following property: if 
$\mathrm{id} = \sigma^{(0)} \rightarrow \sigma^{(1)} \rightarrow \cdots \rightarrow \sigma^{(k)}$ with $0 \leq k \leq \frac{r(r-1)}{2}$, then 
$\sigma^{(i)}=\omega[i]$ for all $0 \leq i \leq k$. In particular, we have $\omega[i-1] \rightarrow \omega[i]$ for all 
$1 \leq i \leq \frac{r(r-1)}{2}$.

\bigskip
 
\begin{lemma} \label{lemma_existence_good_sequence}
Let $i$, $j \in \interv{1}{r}$ be such that $i<j$. Then there exist two integers $k \in \interv{1}{r - 1}$ and 
$e \in \interv{0}{\frac{r(r-1)}{2}-1}$, determined in a unique way by the following properties: $(\omega[e])(k)=i$, $(\omega[e])(k+1)=j$ and 
$\omega[e+1] = \sigma_k \omega[e]$. Namely, we have $k=j-i$ and $e = \frac{(j-1)(j-2)}{2} + (i-1)$.
\end{lemma}

\begin{demo} Left to the reader. \end{demo}

\begin{example} Take $r=6$, $i=2$ and $j=5$. Then we have 
$$_{{{1,2,3,4,5,6} \choose {1,2,3,4,5,6}} \rightarrow {{1,2,3,4,5,6} \choose {2,1,3,4,5,6}} \rightarrow {{1,2,3,4,5,6} \choose {2,3,1,4,5,6}} 
\rightarrow {{1,2,3,4,5,6} \choose {3,2,1,4,5,6}} \rightarrow {{1,2,3,4,5,6} \choose {3,2,4,1,5,6}} \rightarrow{{1,2,3,4,5,6}\choose {3,4,2,1,5,6}}
\rightarrow {{1,2,3,4,5,6} \choose {4,3,2,1,5,6}} \rightarrow {{1,2,3,4,5,6} \choose {4,3,2,5,1,6}} = \omega[e]} ,$$ 
whence $e=7=\frac{(j-1)(j-2)}{2} + (i-1)$ and $k=3=j-i$. \fini
\end{example}							

\begin{prop} \label{prop_existence_good_sequence} Recall Notation \ref{notation_H}. Assume that $(\lambdal,\mul)$ satisfies (\ref{assumption_H}),
$(\gamma,\delta) \neq (0,0)$ and $h \equiv \eta \pmod{nl}$, with 
$\eta \in \set{0,\gamma,\delta,\gamma+\delta}$. Then there exists a good sequence with respect to $(\lambda,\mu)$.
\end{prop}

\begin{demo} 
We construct a good sequence
\begin{equation}
\tag{$*$} v_{\omega.\mu} \stackrel{\bullet}{\rightarrow} \cdots \stackrel{\bullet}{\rightarrow} \textbf{u} 
\stackrel{t}{\rightarrow} \textbf{v} \stackrel{\bullet}{\rightarrow} \cdots \stackrel{\bullet}{\rightarrow} v_{\lambda}
\end{equation} 
as follows.

\begin{itemize}
\item \underline{\emph{Step 1:} construction of 
$v_{\omega.\mu} \stackrel{\bullet}{\rightarrow} \cdots \stackrel{\bullet}{\rightarrow} \textbf{u}$}.
We construct this part of Sequence $(*)$ in terms of the relation $\rightarrow$ on $\Sr$. Let $e$ and $k=x-y$ be the integers given by Lemma 
\ref{lemma_existence_good_sequence} applied with the integers $i:=r+1-x$ and $j:=r+1-y$. Then we 
have the sequence $\omega[0] \rightarrow \cdots \rightarrow \omega[e]$, hence the sequence 
$v_{{\omega[0]}^{-1}.(\omega.\mu)} \stackrel{\bullet}{\rightarrow} \cdots \stackrel{\bullet}{\rightarrow} v_{{\omega[e]}^{-1}.(\omega.\mu)}$ is 
admissible. Put $\textbf{u} := v_{{\omega[e]}^{-1}.(\omega.\mu)} = v_{(\omega{\omega[e])}^{-1}.\mu}$. \\

\item \underline{\emph{Step 2:} construction of $\textbf{u} \stackrel{t}{\rightarrow} \textbf{v}$}.
By assumption on $e$ we have
$$\textbf{u} = v_{m_1} \wedge \cdots \wedge v_{m_{k-1}} \wedge v_{\beta_{x}} \wedge v_{\beta_{y}} \wedge v_{m_{k+2}} \wedge \cdots \wedge v_{m_r},$$
where the $m_i$'s are integers in $B(\mu)$, and the next step of the straightening of $\textbf{u}$
consists in straightening this wedge product with respect to its $k$-th and $(k+1)$-th components, namely $v_{\beta_x} \wedge v_{\beta_y}$. 
Since $(\gamma,\delta) \neq (0,0)$, this elementary straightening involves Rule
$(R_t)$ with $t \in \interv{2}{4}$. Note that by (\ref{eq_1_remark_H_and_beta_numbers}), the wedge product $v_{\beta_y-h} \wedge v_{\beta_x+h}$ is 
ordered. Since $h \equiv \eta \pmod{nl}$, Rule $(R_t)$ shows that this wedge product appears in the linear combination obtained by straightening 
$v_{\beta_x} \wedge v_{\beta_y}$. Put
$$\textbf{v}:=v_{m_1} \wedge \cdots \wedge v_{m_{k-1}} \wedge v_{\beta_y-h} \wedge v_{\beta_x+h} \wedge v_{m_{k+2}} \wedge \cdots \wedge v_{m_r}.$$ 
It is clear that $\textbf{v}$ is obtained from $v_{\lambda}$ by permutation, and the argument above shows that
$\textbf{u} \stackrel{t}{\rightarrow} \textbf{v}$, which completes Step 2. \\

\item \underline{\emph{Step 3:} 
construction of $\textbf{v} \stackrel{\bullet}{\rightarrow} \cdots \stackrel{\bullet}{\rightarrow} v_{\lambda}$}. 
Set $\textbf{v}_1 := \textbf{v}$. If $\textbf{v}_1$ is not ordered, then the elementary straightening 
of $\textbf{v}_1$ gives a linear combination of wedge products, and one of them, say $\textbf{v}_2$, is obtained from $\textbf{v}_1$ by 
permutation. If we apply this device sufficiently many times, we get eventually an ordered wedge product which is of course $v_{\lambda}$. 
This completes Step 3 and the construction of the good sequence $(*)$. \vspace{-8.5mm}
\end{itemize}
\end{demo}

\bigskip

We now prove the converse of Proposition \ref{prop_existence_good_sequence}.

\begin{prop} \label{prop_uniqueness_good_sequence} Recall Notation \ref{notation_H} and assume that $(\lambdal,\mul)$ satisfies 
(\ref{assumption_H}). Assume moreover that there exists a good sequence 
\begin{equation}
v_{\omega.\mu} = v_{\textbf{k}_0} \stackrel{\bullet}{\rightarrow} \cdots \stackrel{\bullet}{\rightarrow} v_{\textbf{k}_e} 
\stackrel{t}{\rightarrow} v_{\textbf{k}_{e+1}} \stackrel{\bullet}{\rightarrow} \cdots \stackrel{\bullet}{\rightarrow} 
v_{\textbf{k}_N} = v_{\lambda}
\end{equation} 
with respect to $(\lambda,\mu)$, with $t \in \interv{2}{4}$. Then this sequence is unique, $t$ is also uniquely determined and moreover we have
$(\gamma,\delta) \neq (0,0)$ and $h \equiv \eta \pmod{nl}$ with $\eta \in \set{0,\gamma,\delta,\gamma+\delta}$.
\end{prop}  

\begin{demo} 
By assumption, $v_{\textbf{k}_e}$ (\resp $v_{\textbf{k}_{e+1}}$) is obtained from $v_{\mu}$ (\resp $v_{\lambda}$) by permutation. 
By (\ref{eq_2_remark_H_and_beta_numbers}), there exist $\sigma \in \Sr$ and $1 \leq k \leq r-1$ such that  
$$v_{\textbf{k}_e} = v_{\sigma.\mu} = v_{\beta_{\sigma^{-1}(1)}} \wedge \cdots \wedge v_{\beta_{\sigma^{-1}(k-1)}} 
\wedge v_{\beta_{x}} \wedge v_{\beta_{y}} \wedge v_{\beta_{\sigma^{-1}(k+2)}} \wedge \cdots \wedge v_{\beta_{\sigma^{-1}(r)}}  ,$$ 
$$v_{\textbf{k}_{e+1}} = v_{\beta_{\sigma^{-1}(1)}} \wedge \cdots \wedge v_{\beta_{\sigma^{-1}(k-1)}} 
\wedge v_{\beta_{y}-h} \wedge v_{\beta_{x}+h} \wedge v_{\beta_{\sigma^{-1}(k+2)}} \wedge \cdots \wedge v_{\beta_{\sigma^{-1}(r)}} \ ,$$
and $v_{\textbf{k}_{e+1}}$ is obtained from $v_{\textbf{k}_e}$ by straightening $v_{\beta_{x}} \wedge v_{\beta_{y}}$ with the rule $(R_t)$. Since
$t \in \interv{2}{4}$, the last conditions of the statement of this proposition hold. It is not hard to see that $k$ and $e$ satisfy the conditions 
of Lemma \ref{lemma_existence_good_sequence} with $i:=r+1-x$ and $j:=r+1-y$, so $k$ and $e$ are determined in a unique way. This determines 
completely the subsequence $v_{\textbf{k}_0} \stackrel{\bullet}{\rightarrow} \cdots \stackrel{\bullet}{\rightarrow} v_{\textbf{k}_e}$. 
Moreover, $t$ is uniquely determined by considering whether $\gamma$ and $\delta$ are zero or not.
The expression of $v_{\textbf{k}_{e+1}}$ given at the beginning of the proof shows that this wedge product is also determined in a unique way.
Let $\tau \in \Sr$ be the unique permutation such that $v_{\textbf{k}_{e+1}} = v_{\tau.\lambda}$. Note that by Condition (i) of Definition 
\ref{def_adjacent_wedges}, there exists at most one admissible sequence $\textbf{V}$ having a given length and starting at a given wedge product 
such that $m(\textbf{V})=0$. As a consequence, the sequence 
$v_{\tau.\lambda} \stackrel{\bullet}{\rightarrow} \cdots \stackrel{\bullet}{\rightarrow} v_{\lambda}$, 
whose length is $\ell(\tau)$, is in turn determined in a unique way.
\end{demo} 

\subsection{Computation of the length modulo $2$ of the good sequence}

We now deal with the technical part of the proof of Theorem \ref{thm_A_J}. The next proposition will be used to show that if
$a'_{\lambdal,\mul}(1)$ and $j^{\prec}_{\lambdal,\mul}$ are nonzero, then both numbers have the same signs. This proposition deals with the only
cases that we have to consider. 

\begin{prop} \label{compatibility_signs_thm_A_J}
Recall Notation \ref{notation_H} and assume that $(\lambdal,\mul)$ satisfies (\ref{assumption_H}). Assume moreover that 
$(\gamma,\delta) \neq (0,0)$ and $h \equiv \eta \pmod{nl}$ with $\eta \in \set{\gamma,\delta}$. Let then \em
\begin{equation}
\textbf{V} = v_{\omega.\mu} \stackrel{\bullet}{\rightarrow} \cdots \stackrel{\bullet}{\rightarrow} v_{\textbf{k}} 
\stackrel{t}{\rightarrow} v_{\textbf{l}} \stackrel{\bullet}{\rightarrow} \cdots \stackrel{\bullet}{\rightarrow} v_{\lambda}
\end{equation}

\noindent \em denote the unique good sequence with respect to $(\lambda,\mu)$ (see Propositions \ref{prop_existence_good_sequence} and 
\ref{prop_uniqueness_good_sequence}). Denote by $N$ the length of this sequence. Then we have

\begin{equation}
(-1)^{N-1} = (-1)^{\haut(\rho)+\haut(\rho')} \varepsilon(\lambdal,\mul) \,\varepsilon,
\end{equation}

\bigskip

\noindent where $\varepsilon(\lambdal,\mul)$ is the sign defined by (\ref{eq_sign_epsilon_lambdal_mul}) and $\varepsilon$ is the sign defined by
\begin{equation} \label{eq_sign_epsilon_compatibility_signs_thm_A_J}
\varepsilon := \left\{ \begin{array}{cl} 1 & \mbox{if} \quad \delta>0 \mbox{ and } h \equiv \delta \pmod{nl}  , \\ -1 & \mbox{otherwise.}
\end{array} \right.
\end{equation}
\end{prop}

\bigskip

\begin{demo} Let $\sigma$, $\tau \in \Sr$ be the permutations defined by $v_{\sigma \omega.\mu} = v_{\textbf{k}}$ and      
$v_{\tau^{-1}.\lambda} = v_{\textbf{l}}$. Then we have $N=\ell(\sigma)+1+\ell(\tau)$, whence $(-1)^{N-1}=\varepsilon(\sigma)\varepsilon(\tau)$.
By Lemma \ref{lemma_existence_good_sequence}, we can compute $\ell(\sigma)$ and then $\varepsilon(\sigma)$; it is however not straightforward to
compute $\ell(\tau)$. We compute only $\varepsilon(\tau)$ by writing $\tau$ as a product of $7$ permutations 
$\sigma^{(1)}, \ldots, \sigma^{(7)}$ whose signs are easily computable.

\begin{itemize}
\item[*] Set first $$\sigma^{(1)}:=\sigma^{-1} \qquad \mbox{and} \qquad \textbf{v}_1:=v_{\sigma^{(1)}.\textbf{l}}\,;$$ 
$\textbf{v}_1$ is thus obtained from $v_{\omega.\mu}$ by replacing $\beta_x$ (located at the $(r+1-x)$-th component)
by $\beta_y-h$, and $\beta_y$ (located at the $(r+1-y)$-th component) by $\beta_x+h$.

\medskip

\item[*] By Statement 1) of Lemma \ref{ribbon_beta_numbers_3}, one of the following cases occurs: 

\begin{itemize} 
\item First case: $d(\beta_x)=d(\beta_y-h)$ and $d(\beta_y)=d(\beta_x+h)$,
\item Second case: $d(\beta_x)=d(\beta_x+h)$, $d(\beta_y)=d(\beta_y-h)$ and $d(\beta_x) \neq d(\beta_y)$.  
\end{itemize}

It is easy to see that the first case occurs if and only if $\delta=0$ or $h \equiv \delta \pmod{nl}$, and that the second case occurs if and only 
if $\delta>0$ and $h \equiv \gamma \pmod{nl}$. Set 

$$\begin{array}{rl}
& \sigma^{(2)}:= \left\{ \begin{array}{cl}
\mathrm{id} & \mbox{if the first case occurs} \\ (r+1-x, r+1-y) & \mbox{if the second case occurs}
\end{array} \right. \\[5mm]
\mbox{and} & \textbf{v}_2:=v_{\sigma^{(2)} \sigma^{(1)}.\textbf{l}} \, ; 
\end{array}$$

in the second case, $\textbf{v}_2$ is obtained from $\textbf{v}_1$ by permuting $\beta_x+h$ and $\beta_y-h$. In either case, $\sigma^{(2)}$ is 
constructed in order to have $\textbf{d}(\sigma^{(2)} \sigma^{(1)}.\textbf{l})=\textbf{d}(\omega.\mu)$ and subsequently
$v(\sigma^{(2)} \sigma^{(1)}.\textbf{l})=v(\omega.\mu)$.     

\medskip

\item[*] Set $$\sigma^{(3)}:=v(\sigma^{(2)} \sigma^{(1)}.\textbf{l})^{-1} = v(\omega.\mu)^{-1} \qquad \mbox{and} \qquad 
\textbf{v}_3:=v_{\sigma^{(3)} \sigma^{(2)} \sigma^{(1)}.\textbf{l}} \, .$$

By remark \ref{remark_permutation_v} applied to $\sigma^{(2)} \sigma^{(1)}.\textbf{l}$, $\textbf{v}_3$ is a wedge product that can be written as 
$$\textbf{v}_3 = v_{\textbf{k}^{(l)}} \wedge \cdots \wedge v_{\textbf{k}^{(1)}}  ,$$ where each $v_{\textbf{k}^{(b)}}$, $1 \leq b \leq l$ is
a wedge product such that each component $v_k$ of $v_{\textbf{k}^{(b)}}$ satisfies $d(k)=b$. In this case and for the rest of the proof, we say 
that $\textbf{v}_3$  is \emph{block-decomposable} and call $v_{\textbf{k}^{(b)}}$ ($1 \leq b \leq l$) the \emph{$b$-th block} of $\textbf{v}_3$.

\medskip

\item[*] Set now $$\sigma^{(4)}:=\omega(\textbf{k}) = \omega(\textbf{l}) \qquad \mbox{and} \qquad 
\textbf{v}_4:=v_{\sigma^{(4)} \sigma^{(3)} \sigma^{(2)} \sigma^{(1)}.\textbf{l}},$$ 
where $\omega(\textbf{k})$ is defined in Section \ref{section_notation_tau_l} 
(the equality $\omega(\textbf{k}) = \omega(\textbf{l})$ comes from Lemma \ref{ribbon_beta_numbers_3}). For $1 \leq b \leq l$ denote by 
$$r_b:= \sharp \set{1 \leq i \leq r \mid d(\alpha_i)=b } = \sharp \set{1 \leq i \leq r \mid d(\beta_i)=b }$$
the number of factors of the block $v_{\textbf{k}^{(b)}}$ (the equality of both numbers defining $r_b$ comes again from Lemma 
\ref{ribbon_beta_numbers_3}). Let $1 \leq b \leq l$. Then $\sigma^{(4)}$ acts on the $b$-th block of $\textbf{v}_3$ as the permutation  
${1, \ldots, r_b} \choose {r_b, \ldots, 1}$. Since $v_{\mu}$ is ordered, we can see that for all $b \in \interv{1}{l} \setminus \set{d,d'}$,
the $b$-th block of $\textbf{v}_4$ is also ordered. 

\medskip
 
\item[*] Write temporarily $\textbf{v}_4 = v_{k_1} \wedge \cdots \wedge v_{k_r}$, and let $i$ (\resp $j$) $\in \interv{1}{r}$ be such that
$k_i=\beta_y-h$ (\resp $k_j=\beta_x+h$). Define $\sigma^{(5)}$ and $\textbf{v}_5$ by
$$\sigma^{(5)}:= \left\{ \begin{array}{cl}
\mathrm{id} & \mbox{if } \delta >0 \\ (i,j) & \mbox{if } \delta=0
\end{array} \right. \qquad \mbox{and} \qquad \textbf{v}_5:=v_{\sigma^{(5)} \sigma^{(4)} \sigma^{(3)} \sigma^{(2)} \sigma^{(1)}.\textbf{l}} \, ;$$
we have $\sigma^{(5)}=\mathrm{id}$ if and only if $\beta_y-h$ and $\beta_x+h$ are in the same block of $\textbf{v}_4$.

\medskip
 
\item[*] Let $\sigma^{(6)} \in \Sr$ be the permutation that acts separately on each block of $\textbf{v}_5$ by reordering it and set 
$$\textbf{v}_6:=v_{\sigma^{(6)} \sigma^{(5)} \sigma^{(4)} \sigma^{(3)} \sigma^{(2)} \sigma^{(1)}.\textbf{l}} \,.$$
Let us describe the action of $\sigma^{(6)}$ more precisely. If $\delta=0$, then $\sigma^{(6)}$ acts on the $d$-th block of $\textbf{v}_5$ as
the permutation $\pi$ from Lemma \ref{ribbon_beta_numbers_2}, and $\sigma^{(6)}$ acts trivially on the other blocks. If $\delta>0$, then 
$\sigma^{(6)}$ acts as the product of two permutations $\sigma'_d$ and $\sigma'_{d'}$, where each $\sigma'_b$, $b \in \set{d,d'}$ acts on the
$b$-th block of $\textbf{v}_5$ as the permutation denoted by $\sigma$ in Lemma \ref{ribbon_beta_numbers_1} and $\sigma'_b$ acts trivially on 
the other blocks. As a consequence, we have in either case $$\varepsilon(\sigma^{(6)})=(-1)^{\haut(\rho)+\haut(\rho')}.$$

\item[*] Finally, put $$\sigma^{(7)}:=v(\lambda) \qquad \mbox{and} 
\qquad \textbf{v}_7:=v_{\sigma^{(7)} \sigma^{(6)} \sigma^{(5)} \sigma^{(4)} \sigma^{(3)} \sigma^{(2)} \sigma^{(1)}.\textbf{l}} \, .$$
Note that $v_{\lambda}$ is ordered, $\textbf{v}_6$ is obtained from $v_{\lambda}$ by permutation, $\textbf{v}_6$ is block-decomposable and all the 
blocks of $\textbf{v}_6$ are ordered. By the remark following the definition of $\sigma^{(3)}$, we have 
$v_{v(\lambda)^{-1}.\lambda}=\textbf{v}_6$, whence $\textbf{v}_7=v_{\lambda}$. 
\end{itemize}

\medskip

As a consequence, we do have $\tau = \sigma^{(7)} \cdots \sigma^{(1)},$ where the $\sigma^{(i)}$'s are defined above, hence 
$(-1)^{N-1}=\varepsilon(\sigma)\prod_{i=1}^{7}{\varepsilon(\sigma^{(i)})}$. By considering different cases we see that
$\varepsilon(\sigma^{(2)})\,\varepsilon(\sigma^{(5)})=\varepsilon$, where $\varepsilon$ is defined by 
(\ref{eq_sign_epsilon_compatibility_signs_thm_A_J}). Moreover, we have 
$$\varepsilon(\sigma^{(4)})=\varepsilon(\omega(\textbf{k}))=\prod_{b=1}^{l}{(-1)^{\frac{r_b(r_b-1)}{2}}} = (-1)^{\kappa(\textbf{d}(\mu))} .$$ 
We then have $\varepsilon(\sigma^{(3)})\,\varepsilon(\sigma^{(4)})\,\varepsilon(\sigma^{(7)})=\varepsilon(\lambdal,\mul)$, whence the result.
\end{demo}

\section{Proof of Theorem \ref{thm_A_J}} \label{proof_thm_A_J}

Let $\lambdal$, $\mul \in \Pilm$, and $\lambda$, $\mu \in \Pi$ be the partitions such that $\lambdal \leftrightarrow \lambda$ and
$\mul \leftrightarrow \mu$. We must show that $a'_{\lambdal,\mul}(1)=2j^{\prec}_{\lambdal,\mul}$. If $\lambdal \not \prec \mul$, then
$a'_{\lambdal,\mul}(1)=0$ by (\ref{eq_barre_unitr}); on the other hand, we have $j^{\prec}_{\lambdal,\mul}=0$ in this case and we are
done. Assume from now on that $\lambdal \prec \mul$. If $(\lambdal,\mul)$ does not satisfy 
(\ref{assumption_H}), then by Remark \ref{rq_CN_existence_bonne_sequence} there cannot exist any good sequence with respect to $(\lambda,\mu)$, so 
by Proposition \ref{express_aprime_lambda_mu} we have $a'_{\lambdal,\mul}(1)=0$ ; on the other hand, by Remark 
\ref{remark_not_H_implies_j_equals_zero} we also have $j^{\prec}_{\lambdal,\mul}=0$ in this case. Assume now that $(\lambdal,\mul)$ 
satisfies (\ref{assumption_H}), and recall Notation \ref{notation_H}. By Proposition \ref{analysis_J1_J2_2}, one of the cases $(J_1)$ or $(J_2)$ 
occurs and Proposition \ref{proposition_J1_J2} then gives the expression of $j^{\prec}_{\lambdal,\mul}$. Moreover, Propositions 
\ref{prop_existence_good_sequence} and \ref{prop_uniqueness_good_sequence} give necessary and 
sufficient conditions on $\gamma$, $\delta$ and $h$ for the existence of a good sequence, in which case it is unique. Proposition 
\ref{express_aprime_lambda_mu} and Rules $(R_2)-(R_4)$ then give the expression of $a'_{\lambdal,\mul}(1)$. In order to compare 
$a'_{\lambdal,\mul}(1)$ and $j^{\prec}_{\lambdal,\mul}$, we have to consider 12 cases depending on the value of $h$ modulo $nl$ and on whether 
$\gamma$ and $\delta$ are zero or not. The results are shown in Figure \ref{fig:3}. Here $N$ is the length of the good sequence if it exists. 
Theorem \ref{thm_A_J} follows by comparing the last two columns of the array and by applying Proposition \ref{compatibility_signs_thm_A_J} if the 
corresponding numbers are nonzero. \cqfd

%
\vspace*{5mm} 
\begin{figure}[htbp]
\begin{center}
\begin{tabular}{|c||c|c|c|}  
\hline
    & \small Number              &  & \\
Case & \small of good & $\ds \frac{a'_{\lambdal,\mul}(1)}{\varepsilon(\lambdal,\mul)}$ & 
	\footnotesize $\ds \frac{j^{\prec}_{\lambdal,\mul}}{(-1)^{\haut(\rho)+\haut(\rho')}}$ \\    
	& \small sequences &   & \\ \hline 
	
$\gamma=\delta=0$ & 0 & 0 & 0 \\ 
\hline 

$\gamma>0$, $\delta=0$, & & & \\
\footnotesize $h \not \equiv 0 \pmod{nl}$, $h \not \equiv \gamma \pmod{nl}$ & 0 & 0 & 0\\
\hline

$\gamma>0$, $\delta=0$, & & \tiny $(-1)^{N-1} \Bigl. \frac{d}{dq} \Bigr|_{q=1} \Bigl ( -(q^{-2}-1)q^{-2i+1} \Bigr)$ & \\
$i:=\frac{h}{nl} \in \mN^{*}$ & 1 & $= 2 (-1)^{N-1}$ & -1\\
\hline

$\gamma>0$, $\delta=0$, & & \tiny $(-1)^{N-1} \Bigl. \frac{d}{dq} \Bigr|_{q=1} \Bigl ( (q^{-2}-1)q^{-2i} \Bigr)$ & \\
$i:=\frac{h-\gamma}{nl} \in \mN$ & 1 & $= -2 (-1)^{N-1}$ & 1\\
\hline

$\gamma=0$, $\delta>0$, & & & \\
\footnotesize $h \not \equiv 0 \pmod{nl}$, $h \not \equiv \delta \pmod{nl}$ & 0 & 0 & 0\\
\hline

$\gamma=0$, $\delta>0$, & & \tiny $(-1)^{N-1} \Bigl. \frac{d}{dq} \Bigr|_{q=1} \Bigl ( -(q^{2}-1)q^{2i-1} \Bigr)$ & \\
$i:=\frac{h}{nl} \in \mN^{*}$ & 1 & $= -2 (-1)^{N-1}$ & 1\\
\hline

$\gamma=0$, $\delta>0$, & & \tiny $(-1)^{N-1} \Bigl. \frac{d}{dq} \Bigr|_{q=1} \Bigl ( (q^{2}-1)q^{2i} \Bigr)$ & \\
$i:=\frac{h-\delta}{nl} \in \mN$ & 1 & $= 2 (-1)^{N-1}$ & 1\\
\hline

$\gamma>0$, $\delta>0$, & & & \\
\scriptsize $h \not \equiv 0 \pmod{nl}$, $h \not \equiv \gamma \pmod{nl}$, & 0 & 0 & 0\\
\scriptsize $h \not \equiv \delta \pmod{nl}$, $h \not \equiv \gamma+\delta \pmod{nl}$ &  &  & \\
\hline

$\gamma>0$, $\delta>0$, & & \tiny $(-1)^{N-1} \Bigl. \frac{d}{dq} \Bigr|_{q=1} \Bigl ( -(q-q^{-1})\frac{q^{2i}-q^{-2i}}{q+q^{-1}} \Bigr)$ & \\
$i:=\frac{h}{nl} \in \mN^{*}$ & 1 & $= 0$ & 0\\
\hline

$\gamma>0$, $\delta>0$, & & \tiny $(-1)^{N-1} \Bigl. \frac{d}{dq} \Bigr|_{q=1} \Bigl ( -(q-q^{-1})\frac{q^{2i+1}+q^{-2i-1}}{q+q^{-1}} \Bigr)$ & \\
$i:=\frac{h-\gamma}{nl} \in \mN$ & 1 & $= -2 (-1)^{N-1}$ & 1\\
\hline

$\gamma>0$, $\delta>0$, & & \tiny $(-1)^{N-1} \Bigl. \frac{d}{dq} \Bigr|_{q=1} \Bigl ( (q-q^{-1})\frac{q^{2i+1}+q^{-2i-1}}{q+q^{-1}} \Bigr)$ & \\
$i:=\frac{h-\delta}{nl} \in \mN$ & 1 & $= 2 (-1)^{N-1}$ & 1\\
\hline

$\gamma>0$, $\delta>0$, & & \tiny $(-1)^{N-1} \Bigl. \frac{d}{dq} \Bigr|_{q=1} \Bigl ( (q-q^{-1})\frac{q^{2i+2}-q^{-2i-2}}{q+q^{-1}} \Bigr)$ & \\
$i:=\frac{h-\gamma-\delta}{nl} \in \mN$ & 1 & $= 0$ & 0\\
\hline 
  	
\end{tabular}
\end{center}
\caption{List of the cases involved in the proof of Theorem \ref{thm_A_J}.}
\label{fig:3}
\end{figure}

\vspace{5mm}
\small \rm Xavier YVONNE, Laboratoire de Math\'ematiques Nicolas Oresme, Universit\'e de Caen, BP 5186, 14032 Caen Cedex, France. \\

\indent \it E-mail address: \tt xyvonne@math.unicaen.fr

\clearpage
\begin{thebibliography}{ABCD}
\bibitem[A1]{A1} \scshape S. Ariki, \rm \emph{On the semisimplicity of the Hecke algebra of $(\mZ/r\mZ) \wr {\mathfrak{S}_n}$}, J. Algebra 
\textbf{169} (1994), 216-225.
\bibitem[A2]{A2} \scshape S. Ariki, \rm \emph{On the decomposition numbers of the Hecke algebra of $G(m,1,n)$}, J. Math. Kyoto Univ. 
\textbf{36} (1996), 789-808. 
\bibitem[AK]{AK} \scshape S. Ariki, K. Koike, \rm \emph{A Hecke algebra of $(\mZ/r\mZ) \wr {\mathfrak{S}_n}$ and construction of its irreducible
representations}, Adv. Math. \textbf{106} (1994), 216-243.
\bibitem[AM]{AM} \scshape S. Ariki, A. Mathas, \rm \emph{The representation type of Hecke algebras of type $B$}, 
Adv. Math. \textbf{181} (2004), 134-159.
\bibitem[BM]{BM} \scshape M. Brou\'e, G. Malle, \rm \emph{Zyklotomische Heckealgebren},
Soci\'et\'e Math\'ematique de France, Ast\'erisque \textbf{212} (1993), 119-189.
\bibitem[DJ]{DJ} \scshape R. Dipper, G. James, \rm \emph{The $q$-Schur algebra}, Proc. London Math. Soc. (3), \textbf{59} (1989), 23-50. 
\bibitem[DJM]{DJM} \scshape R. Dipper, G. James, A. Mathas, \rm \emph{Cyclotomic $q$-Schur algebras}, Math. Z., \textbf{229} (1998), 385-416. 
\bibitem[GJ]{GJ} \scshape O. Gabber, A. Joseph, \rm \emph{Towards the Kazhdan-Lusztig conjecture}, Ann. Sci. \'Ec. Norm. Sup\'er., IV. S\'er. 
\textbf{14} (1981), 261-302.
\bibitem[GGOR]{GGOR} \scshape V. Ginzburg, N. Guay, E. Opdam, R. Rouquier, \rm \emph{On the category $\mathcal{O}$ for rational Cherednik algebras},
Inventiones Math. \textbf{154} (2003), 617-651.
\bibitem[GL]{GL} \scshape J. J. Graham, G. I. Lehrer, \rm \emph{Cellular algebras}, Invent. Math. \textbf{123} (1996), 1-34.
\bibitem[Jac]{Jac} \scshape N. Jacon, \rm \emph{On the parametrization of the simple modules for Ariki-Koike algebras}, J. Math. Kyoto Univ. 
\textbf{44} (2004), 729-767.
\bibitem[JM]{JM} \scshape G. James, A. Mathas, \rm \emph{The Jantzen sum formula for cyclotomic $q$-Schur algebras}, Trans. AMS \textbf{352}, 
(2000), 5381-5404.
\bibitem[Jan]{Jan} \scshape J. C. Jantzen, \rm \emph{Darstellungen halbeinfacher algebraischer Gruppen und zugeordnete kontravariante Formen},
Bonn. Math. Schr., \textbf{67} (1973).
\bibitem[KT]{KT} \scshape M. Kashiwara, T. Tanisaki, \rm \emph{Parabolic Kazhdan-Lusztig polynomials and Schubert varieties}, J. Algebra 
\textbf{249} (2002), 306-325.
\bibitem[Mac]{Mac} \scshape I. G. Macdonald, \rm \emph{Symmetric Functions and Hall Polynomials}, 2nd ed. Oxford University Press (1995).
\bibitem[Mat1]{Mat1} \scshape A. Mathas, \rm \emph{Iwahori-Hecke Algebras and Schur Algebras of the Symmetric Group}, AMS University Lecture Series 
\textbf{15} (1999). 
\bibitem[Mat2]{Mat2} \scshape A. Mathas, \rm \emph{The representation theory of the Ariki-Koike and cyclotomic $q$-Schur algebras}, Adv. Stud. 
Pure Math. \textbf{40} (2004), 261-320.
\bibitem[Ro]{Ro} \scshape R. Rouquier, \rm \emph{$q$-Schur algebras and complex reflection groups, I}, arXiv math.RT/0509252 (2005).
\bibitem[Ry]{Ry} \scshape S. Ryom-Hansen, \rm \emph{The Schaper Formula and the Lascoux, Leclerc and Thibon-algorithm}, 
Lett. Math. Phys. \textbf{64} (2003), 213-219. 
\bibitem[U1]{U1} \scshape D. Uglov, \rm \emph{Canonical bases of higher-level $q$-deformed Fock spaces}, arXiv math.QA/9901032 (1999).
\bibitem[U2]{U2} \scshape D. Uglov, \rm \emph{Canonical bases of higher-level $q$-deformed Fock spaces and Kazhdan-Lusztig polynomials}, 
in Physical Combinatorics ed. M. Kashiwara, T. Miwa, Progress in Math. \textbf{191}, Birkhäuser (2000), arXiv math.QA/9905196 (1999).  
\bibitem[VV]{VV} \scshape M. Varagnolo, E. Vasserot, \rm \emph{On the decomposition matrices of the quantized Schur algebra}, Duke Math. J. 
\textbf{100} (1999), 267-297. 
\bibitem[Y]{Y} \scshape X. Yvonne, \rm \emph{Bases canoniques d'espaces de Fock de niveau sup\'erieur}, Thèse de l'Universit\'e de Caen (2005).

\end{thebibliography}
\end{document}